\documentclass[12pt]{amsart}
\usepackage{amsmath, amssymb, latexsym}

\newcommand{\arx}[1]{\texttt{http://arxiv.org/abs/#1}}
\newcommand{\op}{\operatorname}
\newcommand{\maxfin}{\op{maxfin}}

\newcommand{\Cal}{\mathcal}

\newcommand{\A}{{\Cal A}}
\newcommand{\B}{{\Cal B}}
\newcommand{\BG}{\B_\Gamma}
\newcommand{\BO}{\B_\Omega}

\newcommand{\F}{{\Cal F}}
\newcommand{\J}{{\Cal J}}

\newcommand{\M}{{\Cal M}}
\newcommand{\N}{\naturals}
\newcommand{\NN}{{{}^{\naturals}\naturals}}
\renewcommand{\inf}{[\naturals]^{\infty}}
\renewcommand{\O}{\Cal O}
\renewcommand{\P}{\mathbb P}
\newcommand{\Q}{\rationals}
\newcommand{\R}{\reals}

\newcommand{\UU}{{\Cal U}}

\newcommand{\V}{{\Cal V}}
\newcommand{\Z}{{\mathbb Z}}

\long\def\?#1\QED{}
\long\def\forget#1\forgotten{}

\renewcommand{\b}{{\mathfrak b}}
\renewcommand{\c}{{\mathfrak c}}
\renewcommand{\d}{{\mathfrak d}}

\newcommand{\g}{\gamma}
\renewcommand{\i}{\item}
\renewcommand{\k}{\kappa}

\newcommand{\p}{{\mathfrak p}}
\newcommand{\s}{\sigma}
\newcommand{\w}{\omega}
\newcommand{\x}{\times}

\newcommand{\nin}{\not\in}

\newcommand{\sbst}{\subseteq}
\newcommand{\spst}{\supseteq}
\newcommand{\sm}{\setminus}

\newcommand{\as}{\subseteq^*}
\renewcommand{\pi}{pseudo-intersection}
\renewcommand{\(}{\left(}
\renewcommand{\)}{\right)}

\newcommand{\I}[1]{\emph{#1}}

\newcommand{\cov}{{\sf cov}}
\newcommand{\add}{{\sf add}}
\newcommand{\cof}{{\sf cof}}

\newcommand{\unif}{{\sf non}}
\newcommand{\COV}{{\sf COV}}

\newcommand{\COF}{{\sf COF}}

\long\def\note#1\endnote{\ifNote{\par\medskip\tt Note: #1}\medskip\par
                         \else{}
                         \fi}

\renewcommand{\t}{\tilde}

\newtheorem{thm}{Theorem}
\newtheorem{prop}[thm]{Proposition}
\newtheorem{prob}[thm]{Problem}
\newtheorem{lem}[thm]{Lemma}
\newtheorem{cor}[thm]{Corollary}

\theoremstyle{remark}
\newtheorem{rem}[thm]{Remark}

\newcommand{\be}{\begin{enumerate}}
\newcommand{\ee}{\end{enumerate}}
\newcommand{\bi}{\begin{itemize}}
\newcommand{\ei}{\end{itemize}}

\newcommand{\sone}{{\sf S}_1}    \newcommand{\sfin}{{\sf S}_{fin}}
\newcommand{\ufin}{{\sf U}_{fin}}

\newcommand{\gone}{{\sf G}_1}    \newcommand{\gfin}{{\sf G}_{fin}}

\newcommand{\naturals}{{\mathbb N}}
\newcommand{\reals}{{\mathbb R}}
\newcommand{\rationals}{{\mathbb Q}}
\newcommand{\integers}{{\mathbb Z}}

\newcommand{\gdelta}{{\sf G}_{\delta}}
\newcommand{\fsigma}{{\sf F}_{\sigma}}
\newcommand{\closure}{\overline}

\begin{document}

\title{The combinatorics of Borel covers}

\author{Marion Scheepers}
\address{Department of Mathematics, Boise State University, Boise, Idaho 83725, USA}
\email{mscheep@micron.net}
\thanks{The first author was supported by NSF grant DMS 9971282.}
\author{Boaz Tsaban}
\address{Department of Mathematics and Computer Science, Bar-Ilan University,
Ramat-Gan 52900, Israel}
\email{tsaban@macs.biu.ac.il}
\urladdr{http://www.cs.biu.ac.il/\~{}tsaban}
\thanks{This paper constitutes a part of the second author's doctoral dissertation at
Bar-Ilan University.}

\begin{abstract}
In this paper we extend previous studies of selection principles for families of open covers of sets of real numbers to also include families of countable Borel covers.
The main results of the paper could be summarized as follows:
\begin{enumerate}
\item
{Some of the classes which were different for open covers are equal for Borel covers -- Section 1;}
\item
{Some Borel classes coincide with classes that have been studied under a different guise by other authors -- Section 4.}
\end{enumerate}
\end{abstract}


\maketitle

\section{Introduction}

Let $X$ be a topological space. Let $\O$ denote the collection of all countable open covers of $X$.
According to \cite{gerlitsnagy} an open cover ${\Cal U}$ of $X$ is said to be an $\omega$-cover if $X$ is not a
member of ${\Cal U}$, but for each finite subset $F$ of $X$ there is a $U\in{\Cal U}$ such that $F\sbst U$.
It is shown in \cite{gerlitsnagy} that every $\omega$-cover of $X$ has a countable subset which is an
$\omega$-cover of $X$ if, and only if, all finite powers of $X$ have the Lindel\"of property. All finite powers
of sets of real numbers have the Lindel\"of property. The symbol $\Omega$ denotes the collection of all
\emph{countable} $\omega$-covers of $X$. According to \cite{jmss} and \cite{comb1} an open cover of $X$
is said to be a $\g$-cover if it is infinite and each element of $X$ is a member of all but finitely many
members of the cover. Since each infinite subset of a $\g$-cover is a $\g$-cover, each
$\g$-cover has a countable subset which is a $\g$-cover. The symbol $\Gamma$ denotes the
collection of all \emph{countable} $\g$-covers of $X$.

Let $\A$ and $\B$ be collections of subsets of $X$.
The following two selection hypotheses have a long history for the case when $\A$ and $\B$ are collections of
   topologically significant subsets of a space. Early instances of these can be found in \cite{hurewicz25} and
   \cite{rothberger38}; many papers since then have studied these selection hypotheses in one form or another.
\begin{itemize}
\item[$\sone(\A,\B)$:]{ For each sequence $(A_n:n\in\N)$ of members of $\A$, there is a sequence
   $(b_n:n\in\N)$ such that for each $n$ $b_n\in A_n$, and $\{b_n:n\in\N\}\in\B$.}
\item[$\sfin(\A,\B)$:]{ For each sequence $(A_n:n\in\N)$ of members of $\A$, there is a sequence
   $(B_n:n\in\N)$ such that each $B_n$ is a finite subset of $A_n$, and $\cup_{n\in\N}B_n\in\B$.}
\end{itemize}
   These selection hypotheses are monotonic in the second variable and anti-monotonic in the first. Moreover, each has a
   naturally associated game:

   In the game $\gone(\A,\B)$ ONE chooses in the $n$-th inning an element $O_n$ of $\A$ and then TWO responds by choosing $T_n\in O_n$. They play an inning per natural number. A play $(O_1, T_1, \dots, O_n,T_n,\dots)$ is won by TWO if $\{T_n:n\in\naturals\}$ is a member of $\B$; otherwise, ONE wins. If ONE does not have a winning strategy in $\gone(\A,\B)$, then $\sone(\A,\B)$ holds. The converse is not always true; when it is true, the game is a powerful tool for studying the combinatorial properties of $\A$ and $\B$.

   The game $\gfin(\A,\B)$ is played similarly. In the $n$-th inning ONE chooses an element $O_n$ of $\A$ and TWO responds with a finite set $T_n\sbst O_n$. A play $O_1, T_1,\dots,O_n,T_n,\dots$ is won by TWO if $\cup_{n\in\naturals}T_n$ is in $\B$; otherwise, ONE wins. As above: If ONE has no winning strategy in $\gfin(\A,\B)$, then $\sfin(\A,\B)$ holds; when the converse is also true the game is a powerful tool for studying $\A$ and $\B$.

A third selection hypothesis, introduced by Hurewicz in \cite{hurewicz25}, is as follows:
\begin{itemize}
\item[$\ufin(\A,\B)$:]{ For each sequence $(A_n:n\in\N)$ of members of $\A$, there is a sequence $(B_n:n\in\N)$ such that for each $n$ $B_n$ is a finite subset of $A_n$, and either $\cup B_n = X$ for all but finitely many $n$, or else $\{\cup B_n:n\in\N\}\setminus\{X\}\in\B$.}
\end{itemize}

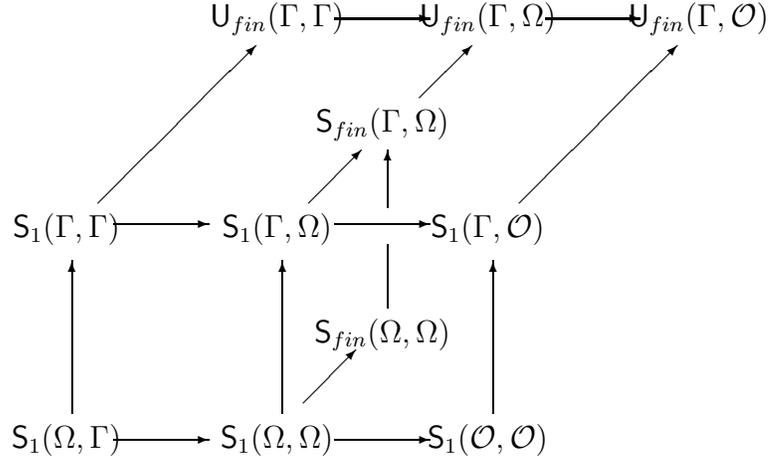
\begin{figure}
\unitlength=.7mm
\begin{picture}(120.00,90.00)(10,10)
\put(20.00,20.00){\makebox(0,0)[cc]
{$\sone(\Omega,\Gamma)$ }}
\put(60.00,20.00){\makebox(0,0)[cc]
{$\sone(\Omega,\Omega)$ }}
\put(100.00,20.00){\makebox(0,0)[cc]
{$\sone(\O,\O)$ }}
\put(20.00,60.00){\makebox(0,0)[cc]
{$\sone(\Gamma,\Gamma)$ }}
\put(60.00,60.00){\makebox(0,0)[cc]
{$\sone(\Gamma,\Omega)$ }}
\put(100.00,60.00){\makebox(0,0)[cc]
{$\sone(\Gamma,\O)$ }}
\put(80.00,40.00){\makebox(0,0)[cc]
{$\sfin(\Omega,\Omega)$ }}
\put(80.00,80.00){\makebox(0,0)[cc]
{$\sfin(\Gamma,\Omega)$ }}
\put(60.00,100.00){\makebox(0,0)[cc]
{$\ufin(\Gamma,\Gamma)$ }}
\put(100.00,100.00){\makebox(0,0)[cc]
{$\ufin(\Gamma,\Omega)$ }}
\put(140.00,100.00){\makebox(0,0)[cc]
{$\ufin(\Gamma,\O)$ }}
\put(70.00,100.00){\vector(1,0){18.00}}
\put(110.00,100.00){\vector(1,0){18.00}}
\put(86.00,85.00){\vector(1,1){10.00}}
\put(25.00,65.00){\vector(1,1){30.00}}
\put(105.00,65.00){\vector(1,1){30.00}}
\put(65.00,65.00){\vector(1,1){10.00}}
\put(28.00,61.00){\vector(1,0){18.00}}
\put(70.00,61.00){\vector(1,0){18.00}}
\put(80.00,45.00){\line(0,1){12.00}}
\put(80.00,64.00){\vector(0,1){11.00}}
\put(20.00,25.00){\vector(0,1){29.00}}
\put(60.00,25.00){\vector(0,1){29.00}}
\put(64.00,27.00){\vector(1,1){10.00}}
\put(100.00,25.00){\vector(0,1){29.00}}
\put(28.00,20.00){\vector(1,0){18.00}}
\put(70.00,20.00){\vector(1,0){18.00}}
\end{picture}
\caption{The open covers diagram \label{basicdiagr}}
\end{figure}

The three classes of open covers above are related: $\Gamma\sbst\Omega\sbst\O$. This and the
properties of the selection hypotheses lead to a complicated diagram depicting how the classes defined this way interrelate. However, only a few of these classes are really distinct, as was shown in \cite{jmss} and \cite{comb1}. Figure \ref{basicdiagr} (borrowed from \cite{jmss}) contains the distinct ones among these classes (it is not known if the class $\sfin(\Gamma,\Omega)$ is $\ufin(\Gamma,\Omega)$, or if it contains $\ufin(\Gamma,\Gamma)$). In this diagram, as in the  ones to follow, an arrow denotes implication.

    Now we consider the following covers of $X$. The symbol $\B$ denotes the family of all \emph{countable} covers of $X$ by \emph{Borel sets}; call elements of $\B$ countable Borel covers of $X$. A countable Borel cover of $X$ is said to be a \emph{Borel $\omega$-cover} of $X$ if $X$ is not a member of it but for each finite subset of $X$ there is a member of the cover which contains the finite set. The symbol $\BO$ denotes the collection of Borel $\omega$-covers of $X$. A countable Borel cover of $X$ is said to be a \emph{Borel $\g$-cover} of $X$ if it is infinite and each element of $X$ belongs to all but finitely many members of the cover. The symbol $\BG$ denotes the collection of Borel $\g$-covers of $X$.
    It is evident that the following inclusions hold:
\[\BG\sbst \BO\sbst \B; \, \Gamma\sbst \BG; \, \Omega\sbst    \BO \mathrm{\ and\ } \O\sbst\B.
\]

    On account of these inclusions and monotonicity properties of the selection principles we have: $\sone(\B,\B)\sbst \sone(\O,\O);
\, \sfin(\B,\B)\sbst \sfin(\O,\O); \, \ufin(\BG,\BG) \sbst \ufin(\Gamma,\Gamma);
\, \sone(\BO,\BG) \sbst \sone(\Omega,\Gamma);$
and so on.

The methods of \cite{jmss} and \cite{comb1} can be used to show that a diagram obtained
from Figure \ref{basicdiagr} by substituting all the open classes by their corresponding
Borel versions summarizes all the interrelationships among these.

But there are big differences about what is provable in these two situations.
For example, it has been shown in \cite{jmss} and \cite{wqn} that there
always is an uncountable set of real numbers in the class $\sone(\Gamma,\Gamma)$ and thus in
$\ufin(\Gamma,\Gamma)$.
According to a result of \cite{miller} it is consistent that no uncountable set of
real numbers has property $\ufin(\BG,\BG)$. Thus it is consistent that some of the
classes which provably do not coincide in the open covers diagram, do coincide in the
Borel covers diagram.

It must be checked which, if any, of the classes in the Borel covers diagram
are provably equal; this is our first task.

\section{Characterizations and equivalence of properties}

In this section we give a number of characterizations for
some of the Borel classes above. In particular, we get that some of the new properties are
equivalent, even though their ``open'' versions are not provably equivalent.

\subsection*{The classes $\sone(\BG,\BG)$, $\sfin(\BG,\BG)$, and $\ufin(\BG,\BG)$}

\begin{thm}\label{borelbounded}
For a set $X$ of real numbers, the following are equivalent:
\begin{enumerate}
\item{$X$ has property $\sone(\BG,\BG)$,}
\item{$X$ has property $\sfin(\BG,\BG)$,}
\item{$X$ has property $\ufin(\BG,\BG)$;}
\item{Every Borel image of $X$ in $\NN$ is bounded.}
\end{enumerate}
\end{thm}
\begin{proof}
We must show that $3\Rightarrow 4$ and $4\Rightarrow 1$.

$3\Rightarrow 4$: This is a Theorem of \cite{BarSch}. In short,
note that the collections ${\Cal U}_n=\{U^n_m : m\in\N\}$, where
$U^n_m = \{f\in\NN : f(n)<m\}$, are open $\g$-covers of $\NN$.
Assume that $\Psi$ is a Borel function from $X$ to $\NN$. Then the
collections $\B_n = \{\Psi^{-1}[U^n_m] : m\in\N\}$ are in $\BG$
for $X$. For all $n$, the sequence $U^n_m$ is monotonically
increasing with respect to $m$. We may assume that for each $n$
$\B_{n+1}$ refines $\B_n$, so that we can use (1) instead of (3)
to get a sequence $\Psi^{-1}[U^n_{m_n}]\in\B_n$ which is in $\BG$
for $X$. Then the sequence $m_n$ bounds $\Psi[X]$.

$4\Rightarrow 1$:
Assume that $\B_n=\{B^n_k : k\in\N\}$, are in $\BG$ for $X$.
Define a function $\Psi$ from $X$ to $\NN$ so that for each $x$ and $n$:
$$\Psi(x)(n) = \min\{k : \(\forall m\geq k\)\ x\in B^n_m\}.$$
Then $\Psi$ is a Borel map, and so $\Psi[X]$ is bounded, say by the sequence $m_n$.
Then the sequence $(B^n_{m_n} : n\in\N)$ is in $\BG$ for $X$.
\end{proof}

\begin{cor}\label{hure1}
For a set $X$ of real numbers, the following are equivalent:
\be
\i $X$ has property $\ufin(\BG,\BG)$.
\i Every Borel image of $X$ has property $\ufin(\Gamma,\Gamma)$.
\ee
\end{cor}
\begin{proof}
An old Theorem of Hurewicz \cite{hurewicz27} asserts
that $X$ has property $\ufin(\Gamma,\Gamma)$ if, and only if,
every continuous image of $X$ in $\NN$ is bounded.
\end{proof}

\begin{thm}\label{bgammadualtob}
For a set $X$ of real numbers the following are equivalent:
\begin{enumerate}
\item{$X$ has property $\sone(\BG,\BG)$.}
\item{Each subset of $X$ has property $\sone(\BG,\BG)$.}
\item{For each measure zero set $N$ of real numbers, $X\cap N$ has property $\sone(\BG,\BG)$.}
\end{enumerate}
\end{thm}
\begin{proof}
$1\Rightarrow 2$:  This follows immediately from Theorem \ref{borelbounded} and the fact that for sets of
real numbers a function on a subspace which is Borel on the subspace, extends to one which is Borel on the whole space.

$3\Rightarrow 1$: Let $X$ be as in $3$, and let $\Psi$ be a Borel function from $X$ to
$\NN$.
We may assume that $X$ is a subset of $[0,1]$, the unit interval (as was shown in \cite{wqn}, the property
$\sone(\Gamma,\Gamma)$ is preserved by countable unions). Let $\Phi$ be a Borel function from $[0,1]$ to
$\NN$ whose restriction to $X$ is $\Psi$.

By Lusin's Theorem choose for each $n$ a closed subset $C_n$ of the unit interval such that $\mu(C_n)\ge 1 - (\frac{1}{2})^n$,
and such that $\Phi$ is continuous on $C_n$. Since $C_n$ is compact, the image of $\Phi$ on $C_n$ is bounded in
$\NN$, say by $h_n$. The set $N=$ [0,1]$\setminus \cup_{n\in\naturals} C_n$ has measure zero, and so
$X\cap N$ has property $\sone(\BG,\BG)$. It follows that the image under $\Psi$ of $X\cap N$ is bounded, say by $h$. Now let
$f$ be a function which eventually dominates each $h_n$, and $h$. Then $f$ eventually dominates each member of $\Psi[X]$.

Since $\Psi$ was an arbitrary Borel function from $X$ to $\NN$, Theorem \ref{borelbounded} implies that $X$
has property $\sone(\BG,\BG)$.
\end{proof}

\begin{prop}\label{bgammaandsigma}
If a set $X$ of real numbers has the $\sone(\BG,\BG)$ property, then it is a $\sigma$-set.
\end{prop}
\begin{proof}
We show that each ${\sf G}_{\delta}$-subset of $X$ is an $\fsigma$-subset. Thus, let $A$ be a $\gdelta$-subset of $X$, say $A = \cap_{n\in\naturals} U_n$ where for all $n$ $U_n\supseteq U_{n+1}$ are open subsets of $X$. Since $X$ is metrizable, each $U_n$ is an $\fsigma$-set. Write, for each $n$,
\[U_n = \cup_{k\in\naturals}C^n_k
\]
where for all $m$, $C^n_m\sbst C^n_{m+1}$ are closed sets. Then for each $n$ $\B_n:=(C^n_m:m\in\naturals)$ is in $\BG$ for $A$. Since $\sone(\BG,\BG)$ is hereditary, $A$ has this property and we find for each $n$ an $m_n$ such that $(C^n_{m_n}:n\in\naturals)$ is a $\g$-cover of $A$. For each $k$ define
\[ F_k:= \cap_{n\ge k} C^n_{m_n}.
\]
Then each $F_k$ is closed and $A=\cup_{k\in\naturals}F_k$.
\end{proof}

According to Besicovitch \cite{besicovitch} a set $X$ of real numbers is \emph{concentrated on} a set $Q$ if for
every open set $U$ containing $Q$, the set $X\setminus U$ is countable.

\begin{cor}\label{sonebgambgamandconcentr}
If an uncountable set of real numbers is concentrated on a countable subset of itself,
then it does not have property $\sone(\BG,\BG)$.
\end{cor}

\subsection*{The classes $\sone(\BG, \B)$, $\sfin(\BG, \B)$, and $\ufin(\BG, \B)$}

\begin{thm}\label{borelhurewicz}
The following are equivalent:
\be
\item{$X$ has property $\sone(\BG,\B)$.}
\item{$X$ has property $\sfin(\BG,\B)$.}
\item{$X$ has property $\ufin(\BG,\B)$.}
\item{No Borel image of $X$ in $\NN$ is dominating.}
\ee
\end{thm}
\begin{proof}
The proof is similar to that of Theorem \ref{borelbounded}.

$3\Rightarrow 4$: Given a Borel function $\Psi$ from $X$ to $\NN$, define $\B_n$ as in
the proof of Theorem \ref{borelbounded}.
Let $A_k$, $k\in\N$, be a partition of $\N$ into infinitely many infinite sets.
From each sequence of covers $\B_n$, $n\in A_k$, we can extract by (1) a cover
$B^n_{m_n}$ $(n\in A_k)$. Taken together, $B^n_{m_n}$ $(n\in\N)$ form a \I{large}
cover of $X$. Recalling that $B^n_{m_n}=\Psi^{-1}[U^n_{m_n}]$, we get that the sequence
$m_n$ witnesses that $\Psi[X]$ is not dominating.

$4\Rightarrow 1$: With notation as in the proof of Theorem \ref{borelbounded}, we get that
if $m_n$ witnesses that $\Psi[X]$ is not dominating, then $(B^n_{m_n} : n\in\N)$ is a (large)
cover of $X$.
\end{proof}

\begin{cor}\label{hure2}
For a set $X$ of real numbers, the following are equivalent:
\be
\i $X$ has property $\ufin(\BG,\B)$.
\i Every Borel image of $X$ in $\NN$ has property $\ufin(\Gamma,\O)$.
\ee
\end{cor}
\begin{proof}
A Theorem of Hurewicz \cite{hurewicz27} asserts that
a set $X$ is $\ufin(\Gamma,\O)$ if, and only if, every continuous image of $X$ in $\NN$ is not dominating.
\end{proof}

\subsection*{The classes $\sone(\BG,\BO)$, $\sfin(\BG,\BO)$, and $\ufin(\BG,\BO)$}

The characterization of these classes is best stated in the language of
filters. Let $\F$ be a filter over $\naturals$. An equivalence
relation $\sim_\F$ is defined on $\NN$ by
$$f\sim_\F g \Leftrightarrow \{n : f(n) = g(n)\} \in \F.$$
The equivalence class of $f$ is denoted $[f]_\F$, and the set
of these equivalence classes is denoted $\NN/\F$.
Using this terminology, $[f]_\F < [g]_\F$ means
$$\{n: f(n) < g(n)\} \in \F.$$

The following combinatorial notion and the accompanying
Lemma \ref{maxfinlemma} will be used to get a technical version of the
filter-based characterization.

For a family $Y\sbst \NN$, define $\maxfin(Y)$ to be the set of
elements $f$ in $\NN$ for which there is a finite set $F\sbst Y$ such that
\[
        f(n) = \max\{h(n):h\in F\}
\]
for all $n$.

\begin{lem}\label{maxfinlemma}
Let $Y\sbst \NN$ be such that for each $n$ the set
$\{h(n):h\in Y\}$ is infinite. Then the following are equivalent:
\begin{enumerate}
\item{$\maxfin(Y)$ is not a dominating family.}
\item{There is a non-principal filter $\F$ on $\naturals$ such that the subset $\{[f]_{\F}:f\in Y\}$ of the reduced product $\NN/\F$ is bounded.}
\end{enumerate}
\end{lem}
\begin{proof}
$1\Rightarrow 2$:
Choose an $h\in \NN$ which is strictly increasing,
and which is not eventually dominated by any element of $\maxfin(Y)$.
For any finite subset $F$ of $Y$, put $f_F(n) = \max\{g(n):g\in F\}$
for each $n$, and then define the set
\[
   A_F = \{n\in\naturals: f_F(n)\le h(n)\}.
\]
Observe that for finite subsets $F$ and $G$ of $Y$, if $F\sbst G$,
then $A_G \sbst A_F$. Thus, the family $\{A_F:F\sbst Y \mbox{ finite}\}$
is a basis for a filter $\F$ on $\naturals$. By the hypothesis on $Y$
this filter is non-principal. It is evident that $[h]/\F$ is an upper
bound for $Y/\F$.\\
$2\Rightarrow 1$: Let $\F$ be a nonprincipal filter on $\naturals$ such
that $Y/\F$ is bounded, and choose a function $h$ in $\NN$
such that for each $f \in Y$ we have $[f]_{\F} < [h]_{\F}$.
Then for each $f\in Y$ the set $\{n: f(n) \le h(n)\}$ is in $\F$ and is
infinite (since $\F$ is non-principal). Since $\F$ has the finite intersection
property it follows that for each finite subset $F$ of $Y$ the set
$S_F = \{n:(\forall f \in F)(f(n)\le h(n))\}$ is in $\F$.
But then $h$ is not eventually dominated by any element of $\maxfin(Y)$.
\end{proof}


\begin{thm}\label{maxfin}
For a set $X$ of real numbers, the following are equivalent:
\be
\i $X$ has property $\sone(\BG,\BO)$,
\i $X$ has property $\sfin(\BG,\BO)$,
\i $X$ has property $\ufin(\BG,\BO)$;
\i For each Borel function $\Psi$ from $X$ to $\NN$, $\maxfin(\Psi[X])$
   is not a dominating family;
\i For each Borel function $\Psi$ from $X$ to $\NN$, either there is a principal filter
   ${\Cal G}$ for which $\Psi[X]/{\Cal G}$ is finite, or else there is a nonprincipal filter ${\Cal F}$ on $\naturals$
   such that the subset $\Psi[X]/{\Cal F}$ of the reduced product $\NN/{\Cal F}$ is bounded.
\ee
\end{thm}
\begin{proof}
$1 \Rightarrow 2 \Rightarrow 3$ are immediate.
We will first show that $3 \Rightarrow 4 \Rightarrow 1$, and then
use Lemma \ref{maxfinlemma} to establish the equivalence of $4$ and $5$.
As in the previous proof,
for any finite subset $F$ of $Y$, put $f_F(n) = \max\{g(n):g\in F\}$
for each $n$.

$3\Rightarrow 4$:
Let $Y=\Psi[X]$. By the upcoming Theorem \ref{borelimages}, $Y$ has property $\ufin(\BG,\BO)$.
For each $n$ and each $k$, define $U^n_k:=\{f: f(n) < k\}$; then set ${\Cal U}_n :=\{U^n_k:k\in\naturals\}$.
Each ${\Cal U}_n$ is a $\g$-cover of $\NN$ since for each $n$ and for $k< j$ we have
$U^n_k\sbst U^n_j$.
Let $A_k$, $k\in\N$, be a partition of $\N$ into infinitely many infinite sets.
From each sequence of $\g$-covers ${\Cal U}_n$, $n\in A_k$, we can use the $\ufin(\BG,\BO)$ property of $Y$ to extract
an $\w$-cover $(U^n_{m_n} : n\in A_k)$.
Then for each finite $F\sbst X$, we have for each $k\in\N$ an $n\in A_k$ such that $\Psi[F]\sbst U^n_{m_n}$, i.e.,
$f_{\Psi[F]}(n)\le m_n$. Thus, the sequence $m_n$ witnesses that $\maxfin(\Psi[X])$
is not a dominating family.

$4 \Rightarrow 1$: Assume that $\B_n=\{B^n_m : m\in\naturals\}$ are in $\BG$ for $X$.
Define a Borel function $\Psi$ from $X$ to $\NN$ so that for each $x$ and $n$:
$$\Psi(x)(n) = \min\{k : \(\forall m\geq k\)\ x\in B^n_m\}.$$
Note that if $F\sbst X$ is finite, then for all $m\geq f_{\Psi[F]}(n)$, $F\sbst B^n_m$.
Let the sequence $m_n$ witness
that $\maxfin(\Psi[X])$ is not dominating. Then for all finite
$F\sbst X$, $F\sbst B^n_{m_n}$ infinitely many times. That is, $(B^n_{m_n} : n\in\naturals)$
is in $\BO$ for $X$.

$4\Rightarrow 5$: There are two cases to consider:\\
{\sf Case 1:} There is an $n$ such that $\{\Psi(x)(n):x\in X\}$ is finite. Then the principal filter generated by $\{n\}$ does the job.\\
{\sf Case 2:} For each $n$ the set $\{\Psi(x)(n):x\in X\}$ is infinite. Apply Lemma \ref{maxfinlemma}.

$5\Rightarrow 4$: Again consider two cases, and apply Lemma \ref{maxfinlemma}.
\end{proof}

\begin{rem}\label{openmaxfin}
{\rm
The implications $1\Rightarrow 2\Rightarrow 3\Rightarrow 4$ and $4\Leftrightarrow 5$
in Theorem \ref{maxfin} can be proved for the
open version of these properties in a similar manner.
The implication $3\Rightarrow 2$ in the open case is counter-exampled
by the Cantor set \cite{jmss}.
We do not know whether the open version of $4\Rightarrow 3$ is true.
}
\end{rem}

This gives the following characterization of $\mathfrak{d}$:
\begin{cor}
For an infinite cardinal number $\kappa$ the following are equivalent:
\begin{enumerate}
\item{$\kappa < \mathfrak{d}$;}
\item{For each subset $X$ of $\NN$ of cardinality at most $\kappa$, there is a non-principal filter
${\Cal F}$ on $\naturals$ such that in the reduced product $\NN/{\Cal F}$ the set $X/{\Cal F}$ is bounded.}
\end{enumerate}
\end{cor}
\begin{proof}
By Theorem \ref{maxfin}, 2 implies 1.
To see that 1 implies 2, consider an infinite $\kappa < \mathfrak{d}$ and a subset $X$ of $\NN$ which is of
cardinality $\kappa$. We may assume that $y\in X$ whenever there is an $x\in X$ such that $y$ differs from $x$ in only finitely many points. Then $\maxfin(X)$ also has cardinality $\kappa$. By Lemma \ref{maxfinlemma}
there exists a nonprincipal filter ${\Cal F}$ on $\naturals$ such that $X/{\Cal F}$ is bounded in $\NN/{\Cal F}$.
\end{proof}

\begin{thm}
For a set $X$ of real numbers, the following are equivalent:
\begin{enumerate}
\item{$X$ has property $\sone(\BG,\BO)$;}
\item{For each Borel mapping $\Psi$ of $X$ into $^{\naturals}\integers$ there is a nonprincipal filter ${\Cal F}$
such that the subring generated by $\Psi[X]/{\Cal F}$ in the reduced power $^{\naturals}\integers/{\Cal F}$ is
bounded below and above.}
\end{enumerate}
\end{thm}
\begin{proof}
That 2 implies 1 is proved as before.
Regarding 1 implies 2: It is evident that if we confine attention to the ring $^{\naturals}\integers$ with pointwise operations,
then a subset $Y$ of it would have property $\sone(\BG,\BO)$ if, and only if, there is a nonprincipal filter
${\Cal F}$ such that $Y/{\Cal F}$ is bounded from below and from above in $^{\naturals}\integers$.
Let $g$ be an element of $\NN$ such that $\Psi[X]/{\Cal F}$ is bounded by $[g]$.
Since the set $\{n\cdot g:n\in\integers\}\cup\{g^n:n\in\naturals\}$ is countable, we find a single $h$
such that for all $n$ $h$ eventually dominates each of $n\cdot g$ and $g^n$. But then in the reduced power
$^{\naturals}\integers/{\Cal F}$ the element $[-h]$ is a lower bound and the element $[h]$ is an upper bound for
the ring generated by $\Psi[X]/{\Cal F}$.
\end{proof}

\subsection*{The class $\sone(\B,\B)$}

The classes $\sone(\B,\B)$ and $\sone(\BG,\BG)$ appear to be each other's ``duals''.

\begin{thm}\label{hereditary} For a set $X$ of real numbers, the following are equivalent:
\begin{enumerate}
\item{$X$ has property $\sone(\B,\B)$;}
\item{Every subset of $X$ has property $\sone(\B,\B)$;}
\item{For each meager set $M\sbst\reals$, $X\cap M$ has property $\sone(\B,\B)$.}
\end{enumerate}
\end{thm}
\begin{proof}
We must show that 1 implies 2, and that 3 implies 1.

{$1\Rightarrow 2$:} This is immediate from the equivalence of $\sone(\B,\B)$ with another notion (see
section 5). However, we give a direct proof.

Let $M$ be a subset of $X$, and assume that $X$ has property $\sone(\B,\B)$. For each $n$ let ${\Cal U}_n$ be a countable cover of $M$ by Borel subsets of $M$. For each $U\in{\Cal U}_n$ let $B_U$ be a Borel subset of $X$ such that $U = M\cap B_U$. Then $X_n:=\cup\{B_U:U\in{\Cal U}_n\}$ is a Borel subset of $X$ since ${\Cal U}_n$ is countable. In turn, $\t X:=\cap_{n\in\N}X_n$ is a Borel subset of $X$.

For each $n$ let $\t{\Cal U}_n$ be $\{B_U:U\in{\Cal U}_n\}\cup\{X\setminus \t X\}$. Then $({\Cal U}_n:n\in\naturals)$ is a sequence of countable Borel covers of $X$. For each $n$ choose a
$V_n\in\t{\Cal U}_n$ such that $\{V_n:n\in\naturals\}$ is a cover of $X$. For each $n$ for which $V_n \neq \t X$, choose $U_n\in{\Cal U}_n$ such that $V_n= B_{U_n}$; for other values of $n$ let $U_n$ be an arbitrary element of ${\Cal U}_n$. Then $(U_n:n\in\naturals)$ covers $M$.

{$3\Rightarrow 1$:} Let $(\B_n:n\in\naturals)$ be a sequence of countable Borel covers of $X$;
enumerate each $\B_n$ as $(B^n_m:m\in\naturals)$.

Since Borel sets have the property of Baire we may choose for each $B^n_m$ an open set
$O^n_m$ and a meager set $M^n_m$ such that
\[B^n_m = (O^n_m\setminus M^n_m) \cup (M^n_m\setminus O^n_m).
\]

Then $A:= \cup_{m,n \in\N} M^n_m$ is a meager set and so $A\cap X$ has property $\sone(\B,\B)$.
For each $n$ such that $n \bmod 3 = 0$, choose a $B^n_{m_n}\in\B_n$ such that $A\cap X$ is covered by these.

For each $n$, $\O_n$, defined to be $\{O^n_m:m\in\naturals\}$, is an open cover of $X\setminus A$.
Let $Q$ be a countable dense subset of $X\setminus A$, and choose for each $n$ with $n \bmod 3 = 1$ an
$O^n_{m_n}$ such that these cover $Q$.

Then the set $B:=X\setminus \cup\{O^n_{m_n}: n\bmod 3 =1\}$ is meager, and so has property $\sone(\B,\B)$.
For each $n$ such that
$n \bmod 3 = 2$, choose an $O^n_{m_n}\in\O_n$ such that these $O^n_{m_n}$'s cover $B$.

Then the sequence $(B^n_{m_n}:n\in\naturals)$ covers $X$.
\end{proof}

Combining of a result from \cite{BarJu} and
\cite{pawlikowskireclaw} with one from \cite{BarSch} yields the
following characterization:

\begin{thm}\label{borelimagec''}
For a set $X$ of real numbers, the following are equivalent:
\be
\i $X$ has property $\sone(\B,\B)$.
\i Each Borel image of $X$ has the Rothberger property $\sone(\O, \O)$.
\ee
\end{thm}

The selection property $\sone(\O,\O)$ manifests itself in several other interesting ways:
these analogues hold also for $\sone(\B,\B)$.

\begin{thm}\label{sonebbgame} For a set $X$ of real numbers, the following are equivalent:
\begin{enumerate}
\item{$\sone(\B,\B)$ holds.}
\item{ONE has no winning strategy in the game $\gone(\B,\B)$.}
\end{enumerate}
\end{thm}
\begin{proof}
We must show that $1\Rightarrow 2$: Let $F$ be a strategy for ONE of the game $\gone(\B,\B)$. Using it, define the following array of Borel subsets of $X$: First, enumerate $F(\emptyset)$, ONE's first move, as $(U_n:n\in\naturals)$. For each response $U_{n_1}$ by TWO, enumerate ONE's corresponding move $F(U_{n_1})$ as $(U_{n_1,n}:n\in\naturals)$. If TWO responds now with $U_{n_1,n_2}$, enumerate ONE's corresponding move $F(U_{n_1},U_{n_1,n_2})$ as $(U_{n_1,n_2,n}:n\in\naturals)$, and so on.

The family $(U_{\tau}: \tau\in{}^{<\omega}\naturals)$ has the property that for each $\tau$ the set $\{U_{\tau\frown n}:n\in\naturals\}$ is a cover of $X$ by Borel subsets of $X$. Moreover, for each function $f$ in $\NN$, the sequence
\[F(\emptyset), U_{f(1)}, F(U_{f(1)}), U_{f(1),f(2)}, F(U_{f(1)},U_{f(1),f(2)}),\dots
\]
is a play of $\gone(\B,\B)$ during which ONE used the strategy $F$. For each such $f$, define $S_f := \cup_{n\in\N}U_{f(1),\dots,f(n)}$. (Thus, $S_f$ is the set of points covered by TWO during a play coded by $f$. We must show that for some such $f$ we have $S_f = X$.

Define the subset $D$ of $X\times {}\NN$ by
\[D:=\{(x,f): x\not\in S_f\}.
\]
Then $D$ is a Borel subset of $X\times {}\NN$. Moreover, for each
$x\in X$ the set $D_x = \{f:x\not\in S_f\}$ is nowhere dense. (To
see this, let $[(n_1,\dots,n_k)]$ be a basic open subset of $\NN$.
Since $\{U_{n_1,\dots,n_k,m}:m\in\naturals\}$ is a cover of $X$
there is an $n_{k+1}$ with $x\in U_{n_1,\dots,n_k,n_{k+1}}$. But
then $[(n_1,\dots,n_k,n_{k+1})]\cap D_x = \emptyset$.) Now recall
from \cite{BarSch} that as $X$ has property $\sone(\B,\B)$ it
follows that $\NN \neq \cup_{x\in X}D_x$ (see section
\ref{connections}). Let $f$ be a function not in $\cup_{x\in X}
D_x$. Then $X=S_f$, and we have defeated ONE's strategy $F$.
\end{proof}

We next show that $\sone(\B,\B)$ is a Ramsey-theoretic property. First observe:
\begin{lem}\label{sonebbissonebomb}
For a set $X$ of real numbers, the following are equivalent:
\be
\i $X$ has property $\sone(\B,\B)$.
\i $X$ has property $\sone(\BO,\B)$.
\ee
\end{lem}
\begin{proof}
The proof for this is like that of Theorem 17 of \cite{comb1}.
\end{proof}

The virtue of $\BO$ for Ramsey-theoretic purposes is that if ${\Cal U}$ is a member of $\BO$,
and if it is partitioned into finitely many pieces, then at least one of these pieces is a member of $\BO$.
This statement is denoted by the abbreviation:
\begin{displaymath}
\mathrm{\ for\ each\ } k\mathrm{,\ } \BO\rightarrow(\BO)^1_k
\end{displaymath}
This is a special case of the more general notation
\[\mathrm{\ for\ all\ }n\mathrm{\ and\ }k\mathrm{\ } \A\rightarrow({\Cal C})^n_k,
\]
which denotes the statement:
\begin{quote}
  For each $n$ and $k$, for each $A\in\A$, and for each $g:[A]^n\rightarrow\{1,\dots,k\}$, there is a $C\sbst A$ such that $C\in{\Cal C}$ and $g$ is constant on $[C]^n$.
\end{quote}

\begin{thm}\label{sonebbandramsey}
For a set $X$ of real numbers the following are equivalent:
\begin{enumerate}
\item{$X$ has property $\sone(\B,\B)$.}
\item{$X$ has the property that for all $k$, $\BO\rightarrow(\B)^2_k$.}
\end{enumerate}
\end{thm}
\begin{proof}
The proof of this is like that of Theorem 4 of \cite{sramsey}.
\end{proof}

\subsection*{The class $\sone(\BO,\BO)$}

It is evident that unions of countably many spaces, each having property $\sone(\B,\B)$,
have property $\sone(\B,\B)$.
\begin{thm}\label{powerssonebbmeanssonebombom}
If all finite powers of $X$ have property $\sone(\B,\B)$, then
$X$ has property $\sone(\BO,\BO)$.
\end{thm}
\begin{proof}
The proof of this is a minor variation on the proof of $(2)\Rightarrow (1)$ of Theorem 3.9 of \cite{jmss}.
\end{proof}

\begin{prob}\label{sonebombompowersareinbb}
Is it true that if $X$ has property $\sone(\BO,\BO)$, then it has property $\sone(\B,\B)$ in all finite powers?
\end{prob}

\subsection*{The class $\sfin(\BO,\BO)$}

It is evident that unions of countably many spaces, each having property $\sone(\BG,\B)$,
have property $\sone(\BG,\B)$.

\begin{thm}\label{powerssonebgbmeanssfinbombom}
If all finite powers of $X$ have property $\sone(\BG,\B)$, then
$X$ has property $\sfin(\BO,\BO)$.
\end{thm}
\begin{proof}
Let $Y=\sum_{k\in\naturals} X^k$. Then by the assumption, $Y$ has property $\sone(\BG,\B)$.
Assume that $\B_n=\{B^n_m : m\in\naturals\}$ are in $\BO$ for $X$.
Define a Borel function $\Psi$ from $Y$ to $\NN$ so that for all $k$, $x_0,\dots,x_{k-1}\in X$, and $n$:
$$\Psi(x_0,\dots,x_{k-1})(n) = \min\{k : \(\forall m\geq k\)\ x_0,\dots,x_{k-1}\in B^n_m\}.$$
By Theorem \ref{borelhurewicz}, the image of $Y$ under $\Psi$ is not dominating. Choose a sequence
$m_n$ witnessing this.
For each $n$, set ${\Cal W}_n:=\{B^n_j:j\le m_n\}$.
Then each ${\Cal W}_n$ is finite, and $\cup_{n\in\naturals}{\Cal W}_n$ is in $\BO$ for $X$.
\end{proof}

\begin{prob}\label{sfinbombompowersareinbgb}
Is it true that if $X$ has property $\sfin(\BO,\BO)$, then it has property $\sone(\BG,\B)$ in all finite powers?
\end{prob}

\subsection*{The class $\sone(\BO,\BG)$}

A standard diagonalization trick gives the following.

\begin{lem}\label{omegagammadiagonal}
The following are equivalent:
\be
\i $X$ has property $\sone(\BO,\BG)$.
\i Every Borel $\w$-cover of $X$ contains a $\g$-cover of $X$.
\ee
\end{lem}
\begin{proof}
The proof of this is like that of the corresponding result in \cite{gerlitsnagy}.
\end{proof}

For the next characterization we need some terminology and notation.
For $a, b\sbst\naturals$, $a\as b$ if $a\sm b$ is finite.
Let $\inf$ denote the set of infinite sets of natural numbers.
$X\sbst \inf$ is \I{centered} if every finite
$F\sbst X$ has an infinite intersection. $a\in \inf$ is a \pi{}
of $X$ if for all $b\in X$, $a\as b$.
$X\sbst\inf$ is a \I{power} if it is centered, but has no \pi{}.

Every countable large Borel cover $\UU=\{U_n : n\in\N\}$ of $X$ is associated with a Borel function
$h_\UU: X\to\inf$, defined
by $h_\UU(x) = \{ n : x\in U_n\}$.

\begin{lem}[{\cite{tsaban}}] Assume that $\UU$ is a cover of $X$. Then:
\be
\i $\UU$ is an $\w$-cover of $X$ if, and only if, $h_\UU[X]$ is centered;
\i $\UU$ contains a $\g$-cover of $X$ if, and only if, $h_\UU[X]$ has a \pi{}.
\ee
\end{lem}

\begin{lem}
The following are equivalent:
\be
\i Every Borel $\w$-cover of $X$ contains a $\g$-cover of $X$.
\i No Borel image of $X$ in $\inf$ is a power.
\ee
\end{lem}
\begin{proof}
$2\Rightarrow 1$: Follows from the preceding lemma.

$1\Rightarrow 2$: Assume that $f:X\to\inf$ is Borel, such that $f[X]$ is centered.
Let $O_n$, $n\in\N$, denote the clopen sets $\{a : n\in a\}$.
As $f[X]$ is centered, $\{O_n : n\in\N\}$ is an $\w$-cover of $f[X]$.
Thus, $\UU=\{f^{-1}(O_n) : n\in\N\}$ is a Borel $\w$-cover of $X$. But $f=h_\UU$,
so we can apply the preceding lemma.
\end{proof}

We thus get the following characterization of $\sone(\BO,\BG)$.
\begin{thm}\label{borelpower}
For a set $X$ of real numbers, the following are equivalent:
\be
\i $X$ has property $\sone(\BO,\BG)$;
\i No Borel image of $X$ in $\inf$ is a power.
\ee
\end{thm}

\begin{cor}\label{omegagammareclaw}
For a set $X$ of real numbers, the following are equivalent:
\be
\i $X$ has property $\sone(\BO,\BG)$.
\i Every continuous image of $X$ has property $\sone(\Omega,\Gamma)$.
\ee
\end{cor}
\begin{proof}
This follows from a Theorem of Rec{\l}aw \cite{irek}, asserting that
$X$ has property $\sone(\Omega,\Gamma)$ if, and only if, no continuous image of $X$ in $\inf$ is a power.
\end{proof}

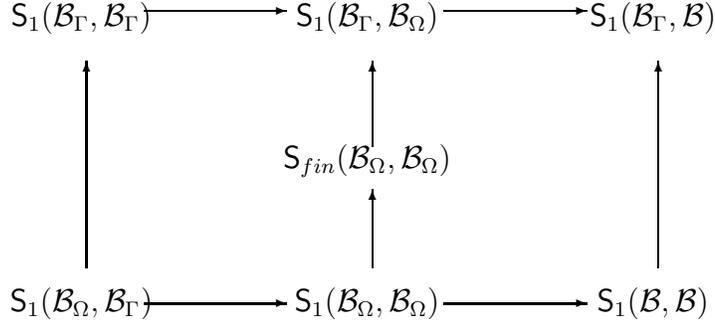
\begin{figure}
\unitlength=.95mm
\begin{picture}(140.00,60.00)(10,10)
\put(20.00,20.00){\makebox(0,0)[cc]
{$\sone(\BO,\BG)$ }}
\put(60.00,20.00){\makebox(0,0)[cc]
{$\sone(\BO,\BO)$ }}
\put(100.00,20.00){\makebox(0,0)[cc]
{$\sone(\B,\B)$ }}
\put(20.00,60.00){\makebox(0,0)[cc]
{$\sone(\BG,\BG)$ }}
\put(60.00,60.00){\makebox(0,0)[cc]
{$\sone(\BG,\BO)$ }}
\put(100.00,60.00){\makebox(0,0)[cc]
{$\sone(\BG,\B)$ }}
\put(60.00,40.00){\makebox(0,0)[cc]
{$\sfin(\BO,\BO)$ }}
\put(28.00,61.00){\vector(1,0){20.00}}
\put(70.00,61.00){\vector(1,0){20.00}}
\put(28.00,20.00){\vector(1,0){20.00}}
\put(70.00,20.00){\vector(1,0){20.00}}

\put(20.00,25.00){\vector(0,1){29.00}}
\put(60.00,25.00){\vector(0,1){11.00}}
\put(60.00,42.00){\vector(0,1){12.00}}
\put(100.00,25.00){\vector(0,1){29.00}}
\end{picture}
\caption{The Surviving Borel Classes \label{boreldiagr}}
\end{figure}

Figure \ref{boreldiagr} summarizes the equivalences proved in this section.

\section{Does Figure \ref{boreldiagr} contain all the provable information about these classes?}
We now consider the question whether we have proved all the equalities that can be proved for these Borel cover classes. It will be seen that the answer is ``Yes"; here is a brief outline of how this follows from the results of the present section:
\begin{enumerate}
\item{According to Corollary \ref{sonebombomandufingammagamma} it is consistent that there is a set of real numbers with property $\sone(\BO,\BO)$, but not property $\sone(\BG,\BG)$. This means that none of the arrows from the left of Figure \ref{boreldiagr} to the middle is reversible.}
\item{According to Theorem \ref{sonebbnotufingamom} it is consistent that there is a set of real numbers in $\sone(\B,\B)$ which is not in $\sone(\BG,\BO)$. This means that none of the arrows from the middle of Figure \ref{boreldiagr} to the right is reversible.}
\item{According to Theorem \ref{sonebgbgandsoneoo} it is consistent that there is a set of real numbers in $\sone(\BG,\BG)$ and not in either of $\sfin(\BO,\BO)$ or $\sone(\B,\B)$. This implies that none of the arrows from the bottom of Figure \ref{boreldiagr} which terminates at the top is reversible.}
\item{According to Theorem \ref{criticalcardinalities} the minimal cardinality of a set of real numbers not having property $\sfin(\BO,\BO)$ is $\mathfrak{d}$, while the minimal cardinality of a set of real numbers not having property $\sone(\B,\B)$ is ${\sf cov}(\mathcal{M})$. Since it is consistent that ${\sf cov}(\mathcal{M}) < \mathfrak{d}$, none of the arrows starting at the bottom row of Figure \ref{boreldiagr} is reversible.}
\end{enumerate}

For a collection $\J$ of 
separable metrizable spaces,
let $\unif(\J)$ denote the minimal cardinality for a
separable metrizable space
which is not a member of $\J$.
We also call $\unif(\J)$ \emph{the critical cardinality} for the class $\J$.

\begin{thm}\label{criticalcardinalities} \
\be
\i $\unif(\sone(\BO,\BG)) = \mathfrak{p}$,
\i $\unif(\sone(\BG,\BG)) = \mathfrak{b}$,
\i $\unif(\sfin(\BO,\BO))  = \unif(\sone(\BG,\BO)) = \unif(\sone(\BG,\B)) = \mathfrak{d}$,
\i $\unif(\sone(\BO,\BO))=\unif(\sone(\B,\B))=\cov({\Cal M})$.
\ee
\end{thm}
\begin{proof}
1 and 2 follow from Theorems \ref{borelpower} and \ref{borelbounded}, respectively.
3 follows from Theorems \ref{borelhurewicz} and \ref{powerssonebgbmeanssfinbombom}.

For 4, we need the following lemma.
\begin{lem}
Let $\J,{\Cal S}$ be collections of 
separable metrizable spaces,
such that $X\in\J$ if, and only if, every Borel image of $X$ is in $\Cal S$.
Then $\unif(\J) = \unif({\Cal S})$.
\end{lem}
\begin{proof}
Since $\J\sbst {\Cal S}$, we have $\unif(\J) \le \unif({\Cal S})$.
Now, let $X$ witness $\unif(\J)$. Then there is a Borel function $\Psi$ on $X$ such that
$\Psi [X]\nin {\Cal S}$. As the cardinality of $\Psi [X]$ cannot be greater than the
cardinality of $X$, we get that $\unif(\J) \ge \unif({\Cal S})$.
\end{proof}
Now, it is well known that $\unif(\sone(\O,\O)) = \cov({\Cal M})$.
Therefore, by Theorem \ref{borelimagec''}, $\unif(\sone(\B,\B))=\cov({\Cal M})$.
Thus, by Theorem \ref{powerssonebbmeanssonebombom},
$\unif(\sone(\BO,\BO))=\cov({\Cal M})$ as well.
\end{proof}

Since it is consistent that $\mathfrak{p} < \cov(\M)$, it is consistent that $\sone(\BO,\BG)$ is not equal
to $\sone(\BO,\BO)$. Similarly the consistency of the inequality $\mathfrak{p} < \mathfrak{b}$ implies that
$\sone(\BO,\BG)$ is not provably equal to $\sone(\BG,\BG)$.

It is consistent that $\mathfrak{b} < \cov(\M)$, and so it is consistent that there is a set of real numbers which has property
$\sone(\BO,\BO)$ but which does not have property $\sone(\BG,\BG)$.

Since it is consistent that $\cov(\M) < \mathfrak{d}$, it is also not provable that $\sfin(\BO,\BO)$ is equal to
either of $\sone(\BO,\BO)$ or $\sone(\B,\B)$.


What the cardinality results do not settle is whether $\sone(\BO,\BO)$ provably coincides with
$\sone(\B,\B)$, or whether any of the three classes associated with the cardinal number $\d$ coincides with another.
They also do not give any indication as to what the interrelationships among two classes might be when their
critical cardinals are equal.
To treat these questions we now consider specific examples which could be constructed on the basis of a variety of
axioms which are consistent. All of the axioms that we use have the form of equality between certain well known cardinal invariants.
Readers who are not familiar with this type of axioms may assume the Continuum Hypothesis instead (in this case,
all of the cardinal invariants become equal to $\aleph_1$).

\subsection*{Special elements of $\sone(\B,\B)$}

A set of real numbers is a \emph{Lusin} set if it is uncountable, but its intersection with each meager set of real numbers
is countable. More generally, for a cardinal $\k$ an uncountable set $X\sbst\R$ is said to be a \I{$\k$-Lusin} set if it has
cardinality at least $\k$, but its intersection with each meager set is less than $\k$.
It is evident that
the smaller the value of $\k$, the harder it is for a set to be a $\k$-Lusin set.
Towards the goal of using as weak hypotheses as possible, this means that we would be interested in
$\k$-Lusin sets
for as large a value of
$\k$ that would allow the conclusion we are aiming at. We now work in the group $^\N\integers$ (which
\emph{topologically} is homeomorphic to the set of irrational numbers), and construct from weak axioms
special elements of $\sone(\B,\B)$.


\begin{lem}\label{covmconstruction}
If $\cov(\M) = \cof(\M)$, and if $Y$ is a subset of $^{\N}\integers$
of cardinality at most $\cof(\M)$, then there is a $\cov(\M)$-Lusin set $L\sbst{{}^{\N}\integers}$
such that $Y\sbst L + L$.
\end{lem}
\begin{proof}
Let $\{y_{\alpha}:\alpha<\cov(\M)\}$ enumerate $Y$. Let $\{M_{\alpha}:\alpha<\cov(\M)\}$ enumerate a cofinal family of meager sets, and construct $L$ recursively as follows: At stage $\alpha$ set $X_{\alpha} = \{a_i:i<\alpha\}\cup\{b_i:i<\alpha\}\cup\bigcup_{i<\alpha}M_i$. Then $(y_{\alpha} - X_{\alpha}) \cup X_{\alpha}$ is a union of fewer than $\cov(\M)$ meager sets. Choose an $a_{\alpha}\in {}^{\N}\integers\setminus((y_{\alpha}-X_{\alpha})\cup X_{\alpha})$.
Evidently, $a_{\alpha}\in (y_{\alpha} - {}^{\N}\integers\setminus X_{\alpha})\cap({}^{\N}\integers\setminus X_{\alpha})$. Thus, choose $b_{\alpha}\in {}^{\N}\integers\setminus X_{\alpha}$ for which $y_{\alpha}-b_{\alpha} = a_{\alpha}$. Then we have $y_{\alpha} = a_{\alpha} + b_{\alpha}$.

Finally, set $L = \{a_{\alpha}:\alpha<\cov(\M)\}\cup\{b_{\alpha}:\alpha<\cov(\M)\}$. Then $L$ is a $\cov(\M)$-Lusin set and $L+L\supseteq Y$.
\end{proof}

The next result is used to show that for $\k$ small enough, $\k$-Lusin sets are in $\sone(\B,\B)$.

\begin{cor}\label{lusinaresonebb}
If $X$ is a ${\sf cov}({\Cal M})$-Lusin set, then it has property $\sone(\B,\B)$.
\end{cor}
\begin{proof}
If $M$ is any meager set, then $M\cap X$ has cardinality less than ${\sf cov}({\Cal M})$,
and thus is in $\sone(\B,\B)$. Now apply Theorem \ref{hereditary}.
\end{proof}

The notion of a Lusin set (i.e., an $\aleph_1$-Lusin set in our current notation) was characterized as follows in
\cite{lusin}: For a topological space $X$ let ${\Cal K}$ denote the collection of ${\Cal U}$ such that ${\Cal U}$
is a family of open subsets of $X$, and $X = \cup\{\overline{U}:U\in {\Cal U}\}$. Then $X$ is a Lusin set if,
and only if, it has property $\sone({\Cal K},{\Cal K})$.

Thus we have:
\begin{cor}\label{sonekkissonebb}
If a set of real numbers has property $\sone({\Cal K},{\Cal K})$, then it has property $\sone(\B,\B)$.
\end{cor}

\begin{thm}\label{sonebbnotufingamom}
If $\cov(\M) = \cof(\M)$, then there is a $\cov(\M)$-Lusin set in $\sone(\B,\B)$ which is not in
$\ufin(\Gamma,\Omega)$.
\end{thm}
\begin{proof}
From the cardinality hypothesis and the fact that $\cov(\M)\le \mathfrak{d}\le \cof(\M)$, we see that there is
in $^{\N}\integers$ a dominating family, say $Y$, of cardinality $\cov(\M)$.
Let $L$ be a  $\cov(\M)$-Lusin set as in Lemma \ref{covmconstruction}, such that $L+L\supseteq Y$. As
$\max\{|f(n)|,|g(n)|\}\geq (|f(n)|+|g(n)|)/2$, we see that for the identity mapping $\Psi$,
$\maxfin(\Psi[L])$ is dominating.
Thus, by Remark \ref{openmaxfin}, $L$ does not have property $\ufin(\Gamma,\Omega)$.

By Corollary \ref{lusinaresonebb} $L$ has property $\sone(\B,\B)$.
\end{proof}

This in particular implies that $\sone(\BO,\BO)$ is not provably equivalent to $\sone(\B, \B)$.

\subsection*{Special elements of $\sone(\BO,\BO)$}

Now that we have clarified most of the interrelationships among the Borel classes, we consider how the Borel classes
are related to the classes in Figure \ref{basicdiagr}.
We have just seen that $\sone(\B,\B)$ need not be contained in $\ufin(\Gamma,\Omega)$,
even when the critical cardinalities for sets not belonging to these classes are the same.

Next we treat $\sone(\BO,\BO)$ and $\ufin(\Gamma,\Gamma)$. We show  how to use the
Continuum Hypothesis to construct a Lusin set which has property $\sone(\BO,\BO)$.
Since it is a Lusin set, it does not satisfy $\ufin(\Gamma,\Gamma)$.

In our construction we use the \emph{ad hoc} concept of an $\omega$-\emph{fat} collection of Borel sets. A collection ${\Cal U}$ of Borel sets is said to be \emph{fat} if for each nonempty open interval $J$ and for each dense $\gdelta$-set $G$ there is a $B\in {\Cal U}$ such that $B\cap G\cap J\neq\emptyset$.
It is said to be $\omega$-\emph{fat} if: for each dense $\gdelta$-set $G$ and for every finite family ${\Cal F}$ of nonempty open sets there is a $B\in{\Cal U}$ such that for each $J\in{\Cal F}$, $B\cap J\cap G$ is nonempty.

A number of facts about these $\omega$-fat families of Borel sets will
play a crucial role in our construction. For ease of reference we state
these as lemmas and give proofs where it seems necessary.
\begin{lem}\label{elementary} Let ${\Cal U}$ be an $\omega$-fat family
consisting of countably many Borel sets.
\begin{enumerate}
\item{For each partition of ${\Cal U}$ into two pieces, at least one of
the pieces is $\omega$-fat.}
\item{If ${\Cal U}$ is a Borel $\omega$-cover of the set $X$ and $F$ is a
finite subset of $X$, then $\{U\in{\Cal U}:F\sbst U\}$ is an $\omega$-fat
Borel $\omega$-cover of $X$.}
\end{enumerate}
\end{lem}

\paragraph*{Added in proof}
As stated, item (2) of Lemma \ref{elementary} is wrong: Let ${\Cal
U} = \{\R\sm\Z\} \cup [\Z]^{<\omega}$. Then $\Cal U$ is an
$\omega$-fat $\omega$-cover of $\Z$. But for any nonempty finite
subset $F$ of $\Z$, the collection $\{U \in {\Cal U} : F \subset
U\}$ is not $\omega$-fat. However, if $X$ is a Lusin set such that
for each nonempty basic open set $G$, $X\cap G$ is uncountable,
then (some minor modification of) item (2) of this Lemma holds.
As the special set $X$ which we will construct is a Lusin set, we can easily make sure
that it has the required property and the proof works. This idea
is extended and explained further in \cite[full version]{huremen2}.

\begin{lem}\label{fatborelcontainsdensegdelta} If $\B$ is a countable fat Borel family, then there is a dense $\gdelta$-set contained in $\cup\B$.
\end{lem}
\begin{proof}
Since $B=\cup\B$ is a Borel set, it has the property of Baire. Let $U$ be an open set such that $(U\setminus B) \cup (B\setminus U)$ is meager. Then $U$ is dense, for let $G$ be a dense $\gdelta$ disjoint from that meager set, and let $J$ be a nonempty open interval. Then $J\cap G\cap B$ is nonempty. But $B= (B\setminus U) \cup (B\cap U)$, so that $(B\cap U)\cap J$ is nonempty.

Now $\reals\setminus U$ is nowhere dense, and we may assume that $G$ is also disjoint from this nowhere dense set. But then $G\sbst B$.
\end{proof}

\begin{lem}\label{shrinkinglemma} If ${\Cal U}$ is a countable $\omega$-fat family of Borel sets and
${\Cal F}$ is a finite nonempty family of nonempty open intervals, then there are a $U\in{\Cal U}$ and for each
$J\in{\Cal F}$ a nonempty open interval $I_J\sbst J$ such that the set $U\cap I_J$ is comeager in $I_J$.
\end{lem}
\begin{proof}
Towards proving the contrapositive, take a countable $\omega$-fat family ${\Cal U}$ of Borel sets, and a finite nonempty family ${\Cal F}$ of nonempty open intervals such that:

For each $U\in{\Cal U}$ there is a $J_U\in{\Cal F}$ such that for each nonempty open interval $I\sbst J_U$ the set $U\cap I$ is not comeager in $I$. Fix such a $J_U$ for each $U\in{\Cal U}$.

Since $U\cap J_U$ is a Borel set, it has the property of Baire. Choose an open set $V\sbst J_U$ such that $(V\setminus(U\cap J_U)) \cup ((U\cap J_U)\setminus V)$ is meager. If $V$ is nonempty, then the meagerness of $V\setminus(U\cap J_U)$ implies that $U\cap V$ is comeager in $V$, contradicting the choice of $U$ and $J_U$. Thus, $V$ is empty, and we find that $U\cap J_U$ is meager. Let $G_U$ be a dense $\gdelta$-set disjoint from $U\cap J_U$.

The set $G = \cap_{U\in{\Cal U}} G_U$ is an intersection of countably many dense $\gdelta$-sets, so is a dense $\gdelta$-set. But then $G$ and ${\Cal F}$ witness that ${\Cal U}$ is not $\omega$-fat.
\end{proof}

\begin{lem}\label{countablestuff} Let $S$ be a countably infinite set and let $(F_n:n\in\naturals)$ be an ascending sequence of finite sets with union equal to $S$. If $({\Cal U}_n:n\in\naturals)$ is a sequence of Borel $\omega$-covers of $S$ such that for each $n$ the set $\{U\in {\Cal U}_n: F_n\sbst U\}$ is $\omega$-fat, then there is a sequence $(U_n:n\in\naturals)$ such that for each $n$ $U_n\in{\Cal U}_n$, $\{U_n:n\in\naturals\}$ is a Borel $\g$-cover of $S$, and $\{U_n:n\in\naturals\}$ is $\omega$-fat.
\end{lem}
\begin{proof}
Let $S$, the $F_n$'s, and the ${\Cal U}_n$'s be as in the hypotheses. We may assume for each $n$ that for all $U\in{\Cal U}_n$ we have $F_n\sbst U$. Let $(J_n:n\in\naturals)$ be an enumeration of the nonempty open intervals with rational endpoints.

Consider $n$. Since ${\Cal U}_n$ is $\omega$-fat, choose a $U_n\in{\Cal U}_n$ and
for each $i\le n$ an open nonempty interval $I^i_n\sbst J_i$ such that $I^i_n\cap U_n$ is comeager in $I^i_n$.

Then the sequence $(U_n:n\in\naturals)$ is as desired. To see this, let $G$ be any dense $\gdelta$-set and let $R_1,\,\dots,\, R_n$ be nonempty open intervals. Choose $m$ so large that for each $i\le n$ there is a $j\le m$ with $J_j\sbst R_i$. When we chose $U_m$ it was done so that for some open nonempty intervals $I_j,\,j\le m$ we had $I_j\sbst J_j$ and $U_m\cap I_j$ is comeager in $I_j$, whence $U_m\cap G\cap I_j$ is comeager in $I_j$. But then for each $r\le n$, $U_m\cap G\cap R_r$ is nonempty.
\end{proof}

\begin{lem}\label{ppointlemma}
If $({\Cal U}_n:n\in\naturals)$ is a sequence of countable $\omega$-fat families of Borel sets such that for
each $n$ ${\Cal U}_{n+1}\sbst {\Cal U}_n$, then there is a countable $\omega$-fat family
$\{U_n : n\in\N\}$
of Borel sets such that for each $n$, $U_n\in{\Cal U}_n$.
\end{lem}
\begin{proof}
Let $J_1, J_2, \dots, J_n,\dots$ be a bijective enumeration of a basis for the topology of $\reals$.
Recursively choose for each $n$ sequences $(I^n_k:k\in\naturals)$ of nonempty open intervals, and for each $n$ a
$U_n\in{\Cal U}_n$ such that:
\begin{enumerate}
\item{For $k<n$ we have $I^n_k = J_n$;}
\item{For $k\geq n$ we have $I^n_k\sbst J_n$ and $U_k\cap I^n_k$ is comeager in $I^n_k$.}
\end{enumerate}
This is possible on account of Lemma \ref{shrinkinglemma}. We claim that ${\Cal U}:=\{U_n:n\in\naturals\}$ is $\omega$-fat.

For let $G$ be a dense $\gdelta$-set and let $R_1,\dots, R_k$ be nonempty open intervals. Choose from the basis intervals $J_{n_1},\dots, J_{n_k}$ such that $n_1<\dots<n_k$ and for $1\le i\le k$ we have $J_{n_i}\sbst R_i$. Let $m$ be larger than $n_k$. Then for $1\le i\le k$ we have:
$U_m\cap I^{n_i}_m$ contains a dense $\gdelta$-subset of $I^{n_i}_m$ and so has nonempty intersection with the dense $\gdelta$-set $G$. Since for each $i$ we have $I^{n_i}_m\sbst R_i$ we see that $U\cap R_i\cap G$ is nonempty.
\end{proof}

\begin{lem}\label{targetinglemma} Let $G$ be a dense $\gdelta$ set and let $J$ be a nonempty open interval. If for each $n$ ${\Cal U}_n$ is a countable $\omega$-fat family of Borel sets, then there is an $x\in J\cap G$ such that for each $n$ the set $\{U\in{\Cal U}_n: x\in U\}$ is $\omega$-fat.
\end{lem}
\begin{proof}
For each $n$ let ${\Cal U}_n$ be a countable $\omega$-fat family of Borel sets. Let $J$ be a nonempty open interval, and let $G$ be a dense $\gdelta$-set.

Let $(J_n:n\in\naturals)$ bijectively enumerate a base for the topology of $\reals$, and write $G = \cap_{n\in\naturals}V^1_n$, where $V^1_1\supseteq V^1_2\supseteq\dots$ are dense open sets. Also, write $R_1:= J$. We may assume that the closure of $J$ is compact.

Recursively construct four sequences $((U^i_n:i\le n):n\in\naturals)$, $((I^i_n:i\le n):n\in\naturals)$, $(R_n:n\in\naturals)$ and $((V^i_n:n\in\naturals):i\in\naturals)$, such that the following requirements are satisfied for each $n$:
\begin{enumerate}
\item{For all $k\leq n$, $U^k_n\in{\Cal U}_k\setminus\{U^i_j:i,j< n\}$;}
\item{For each $i\le n$, $I^i_n\sbst J_i$ is a nonempty open interval such that $I^i_n\cap(\cap_{j\le n}U^j_n)$ is comeager in $I^i_n$;}
\item{$R_{n+1}$ is a nonempty open interval with closure contained in $(\cap_{i\leq n}V^i_{n+1})\cap R_n$;}
\item{$R_{n+1}\cap(\cap_{i\leq n}U^i_n)$ is comeager in $R_{n+1}$;}
\item{$V^n_m\supset V^n_{m+1}$ for all $m$ are dense open subsets of $R_n$;}
\item{$R_{n+1}\cap(\cap_{i\leq n}U^i_n)\supseteq \cap_{m\in\naturals}V^{n+1}_m$.}
\end{enumerate}

To see that this recursion can be carried out, first consider $n=1$: Here we already have $R_1$ and each $V^1_n$ specified. Consider $J_1$ and $R_1$, and ${\Cal U}_1$. Apply Lemma \ref{shrinkinglemma} to choose $U^1_1\in{\Cal U}_1$ and intervals $I^1_1$ and $R_2$ such that $\closure{R_2}\sbst R_1\cap V^1_1$ and $U^1_1\cap R_2$ is comeager in $R_2$ and $U^1_1\cap I^1_1$ is comeager in $I^1_1$. Since $U^1_1\cap R_2$ is comeager in $R_2$, choose a descending sequence $(V^2_n:n\in\naturals)$ of open dense subsets of $R_2$ such that $R_2\cap U^1_1\supseteq \cap_{m\in\naturals} V^2_m$. Thus for $n=1$ sets as required by the five recursion specifications have been found.

Suppose now that $n\ge 1$ and that the recursion has been carried through for $n$ steps. Consider $R_n$, $J_1,\,\dots,\, J_n$, and ${\Cal U}_1,\,\dots,\,{\Cal U}_n$.

Choose for $i\le n+1$ sets $U^i_{n+1}\in{\Cal U}_i\setminus\{U^j_k:j,k\le n\}$ and $R_{n+1}$ an open nonempty interval with closure contained in $R_n\cap(\cap_{i\le n}V^i_{n+1})$, as well as open nonempty intervals $I^i_{n+1}$, $i\le n+1$, such that for each $i$ $I^i_{n+1}\sbst J_i$, and $\cap_{k\le n+1}U^k_{n+1}\cap I^i_{n+1}$ is comeager in $I^i_{n+1}$, and $\cap_{k\le n+1} U^k_n\cap R_{n+1}$ is comeager in $R_{n+1}$. This can be done on account of Lemma \ref{shrinkinglemma}. Then let $(V^{n+1}_m:m\in\naturals)$ be a descending sequence of sets open and dense in $R_{n+1}$ such that $R_{n+1}\cap(\cap_{k\le n+1}U^k_{n+1})\supseteq \cap_{m\in\naturals}V^{n+1}_m$.

This shows how to continue the recursion to the next step.

With the recursive procedure completed, for each $n$ put ${\Cal V}_n=\{U^n_k:k\ge n\}$. By the compactness of $\closure{R}_1$, and by specification 3 of the recursion, $\cap_{n\in\naturals}R_n$ is nonempty. Let $x$ be an element of this intersection.

We claim that each ${\Cal V}_n$ is an $\omega$-fat subset of ${\Cal U}_n$, and that for each $V\in{\Cal V}_n$, we have $x\in V\cap J\cap G$.

To see that ${\Cal V}_n$ is $\omega$-fat, let a dense $\gdelta$-set $H$ and a finite set ${\Cal F}$ of nonempty open intervals be given. Choose $m>n$ so large that there is for each $F\in{\Cal F}$ a $J_i$ with $i\le m$ such that $J_i\sbst F$. Then $U^n_m$ was chosen so that for each of the nonempty open intervals $I^i_m\sbst J_i$, we have $U^n_M\cap I^i_m$ comeager in $I^i_m$. But then as $H$ is a comeager set of reals, we have for each $i\le m$ that $U^n_m\cap I^i_m\cap H$ is nonempty. This implies that for each $F\in{\Cal F}$, $U^n_m\cap F\cap H$ is nonempty.

To see that $x$ is a member of each element of ${\Cal V}_n$, consider a $U^n_m\in{\Cal V}_n$. We have $U^n_m\cap R_m\supseteq \cap_{j\in\naturals} V^m_j$. But for each $j\geq m+1$ we have $R_{j+1}\sbst V^m_j$, and as $x$ is in the intersection of the $R_j$'s, it is in the intersection of the $V^m_j$'s, so in $U^n_m$.
\end{proof}

\begin{lem}\label{cocofinalgdeltadense}
If $\add(\M)=\c$, then there exists a family $(G_{\alpha}:\alpha<\aleph_1)$ of dense $\gdelta$-sets
of reals, such that:
\begin{itemize}
\item{For each dense $\gdelta$-set $G$ there is an $\alpha$ with $G_{\alpha}\sbst G$;}
\item{For $\alpha<\beta<\c$ we have $G_{\beta}\sbst G_{\alpha}$.}
\end{itemize}
\end{lem}
\begin{proof}
Let $(M_\alpha : \alpha<\c)$ be a cofinal family of meager sets. We define by induction on $\alpha<\c$
a monotonically increasing sequence $(\t M_\alpha : \alpha<\c)$ of of $\fsigma$ meager sets as follows:
At stage $\alpha$, let $\hat M_\alpha = \cup_{i<\alpha} \t M_i$.
As $\alpha<\add(\M)$, $\hat M_\alpha$ is meager, so let $\t M_\alpha$ be an $\fsigma$ meager
set containing $\hat M_\alpha$.

By the Baire category Theorem, complements of meager sets in $\reals$ are dense. Thus,
setting for each $\alpha$ $G_\alpha = \reals \setminus \t M_\alpha$ yields the
desired sequence.
\end{proof}

\begin{thm}[CH]\label{existencesonebombom}
There is a $\c$-Lusin set which has property $\sone(\BO,\BO)$.
\end{thm}
\begin{proof}
Let $(G_{\alpha}:\alpha<\c)$ be as in Lemma \ref{cocofinalgdeltadense}.
Let $(({\Cal U}^{\alpha}_n:n\in\naturals):\alpha<\c)$ list all $\omega$-sequences where each term is an $\omega$-fat countable family of Borel sets. We shall now recursively construct the desired Lusin set $X$ by choosing for each $\alpha$ a countable dense set $X_{\alpha}$ to satisfy certain requirements, and then setting $X = \cup_{\alpha<\c} X_{\alpha}\cup\rationals$. Together with each $X_{\alpha}$ we shall choose a sequence $(U^{\alpha}_n:n\in\naturals)$ of Borel sets and a sequence $(S_\g(\alpha):\g<\c)$ of infinite subsets of $\naturals$ such that:
\begin{enumerate}
\item{Whenever $\g<\beta<\c$, then $S_\g(\beta)=\naturals$;}
\item{For each $\beta<\c$, for $\g<\nu<\c$ we have
$S_{\nu}(\beta)\sbst^* S_\g(\beta)$;}
\item{For all $\beta$ and $\g$, $\{U^{\beta}_n:n\in S_\g(\beta)\}$ is an $\omega$-fat
$\g$-cover of $\rationals\cup(\cup_{\nu\leq\g}X_{\nu})$.}
\item{For any $\alpha$, if some ${\Cal U}^{\alpha}_n$ is not an $\omega$-cover of $\rationals\cup(\cup_{\nu<\alpha}X_{\nu})$, then for each $n$ we have $U^{\alpha}_n = \reals$;}
\item{If for each $n$ ${\Cal U}^{\alpha}_n$ is an $\omega$-cover of $\rationals\cup(\cup_{\nu<\alpha}X_{\nu})$, then for each $n$ we have $U^{\alpha}_n\in{\Cal U}^{\alpha}_n$, and $\{U^{\alpha}_n:n\in\naturals\}$ is an $\omega$-fat $\g$-cover of $\rationals\cup(\cup_{\nu<\alpha}X_{\nu})$;}
\item{For each $\alpha$, $X_{\alpha}\sbst G_{\alpha} \setminus (\rationals \cup (\cup_{\nu<\alpha}X_{\nu}))$  is dense in $\reals$.}
\end{enumerate}

Before showing that this can be accomplished, we show that constructing $X$ to satisfy these requirements is sufficient.
Thus, let $X$ be obtained like this. Let $({\Cal U}_n:n\in\naturals)$ be a sequence of countable Borel
$\omega$-covers of $X$. Since each $X_{\alpha}$ is dense and contained in $G_{\alpha}$ it follows that for each
$n$ ${\Cal U}_n$ is $\omega$-fat. Thus, for some $\beta$ we have
$({\Cal U}_n:n\in\naturals) = ({\Cal U}^{\beta}_n:n\in\naturals)$. Since each ${\Cal U}^{\beta}_n$ is
an $\omega$-cover of $X$, it is an $\omega$-cover of $\rationals\cup(\cup_{\g<\beta}X_\g)$,
and thus is as in 5. Let $F$ be a finite subset of $X$ and choose a $\beta>\alpha$ such that
$F\sbst \rationals\cup(\cup_{\g\leq\beta} X_\g)$. By 3
$\{U^{\alpha}_n:n\in S_{\beta}(\alpha)\}$ is a $\g$-cover of
$\rationals\cup(\cup_{\g\leq\beta} X_\g)$, whence for some $n$
$F\sbst U^{\alpha}_n$. It follows that $\{U^{\alpha}_n:n\in\naturals\}$ is an $\omega$-cover of
$X$, as desired.

Now the recursive construction: Fix $\rationals$, the set of rational numbers, and ask: Is $({\Cal U}^0_n:n\in\naturals)$ a sequence of $\omega$-covers of $\rationals$?
\begin{itemize}
\item[No:]{Then for each $n$ set $U^0_n = \reals$, choose $X_0\sbst G_0\setminus \rationals$ countable and dense, and put $S_0(0) = \naturals$.}
\item[Yes:]{For each $n$ choose a $U^0_n\in{\Cal U}^0_n$ such that $\{U^0_n:n\in\naturals)$ is an $\omega$-fat $\g$-cover of $\rationals$. Repeatedly apply Lemma \ref{targetinglemma} to recursively choose numbers $x_1\in J_1\cap G_0\setminus\rationals$ and $x_{n+1}\in J_{n+1}\cap G_0\setminus(\rationals\cup\{x_1,\dots,x_n\})$ such that:
${\Cal V}_1:=\{U^0_n:x\in U^0_n\}$ is an $\omega$-fat family of Borel sets, and for each $n$ ${\Cal V}_{n+1}:=\{U^0_m\in{\Cal V}_n:x_{n+1}\in U^0_m\}$ is an $\omega$-fat family of Borel sets. In the end put $X_0 = \{x_n:n\in\naturals\}$, and choose by Lemma \ref{ppointlemma} a ${\Cal V}\sbst {\Cal V}_1$ such that ${\Cal V}$ is $\omega$-fat, and for each $n$ also ${\Cal V}\sbst^*{\Cal V}_n$. Finally set $S_0(0) = \{n: U^0_n\in{\Cal V}\}$. Observe that $\{U^0_n:n\in S_0(0)\}$ is a $\g$-cover of $\rationals\cup X_0$.}
\end{itemize}

This shows that the six recursive requirements are satisfiable for $\alpha=0$. Assume now that $\alpha>0$ is given,
and for each $\beta<\alpha$ we already have $X_{\beta}$ as well as the sequence $(U^{\beta}_n:n\in\naturals)$
and $(S_\g(\beta):\g<\alpha)$ such that the six recursive requirements are satisfied. To verify that stage
$\alpha$ can then be carried out, do the following. First, for all $\beta<\alpha$ define
$S_{\beta}(\alpha) = \naturals$. Also, using Lemma \ref{ppointlemma}, choose for each $\beta<\alpha$ an infinite
set $S_{\beta}\sbst\naturals$ such that for all $\g<\alpha$ we have
$S_{\beta}\sbst^* S_\g(\beta)$,
and such that $\{U^{\beta}_n:n\in S_{\beta}\}$ is an
$\omega$-fat $\g$-cover
\forget
\\
 ({\huge OOPS! It seems to me that we are using Lemma 2.10 here, and it works for COUNTABLE, while the union following here may be uncountable ... should we go back to CH, or can this be fixed?? Same with hypothesis of Corollary 2.15})\\
\forgotten
 of $\cup_{\g<\alpha}X_\g\cup\rationals$.

Consider $({\Cal U}^{\alpha}_n:n\in\naturals)$ and ask: Is each ${\Cal U}^{\alpha}_n$ an $\omega$-cover of $\cup_{\g<\alpha}X_\g\cup\rationals$?
\begin{itemize}
\item[No:]{ Then for each $n$ put $U^{\alpha}_n=\reals$, and declare $S_{\alpha}(\alpha) = \naturals$.
Next we choose $X_{\alpha}$ recursively as follows from
$H_{\alpha}:= G_{\alpha}\setminus(\cup_{\g<\alpha}X_\g\cup\rationals)$:
By Lemma \ref{targetinglemma} choose an $x_1\in J\cap H_{\alpha}$ such that for each $\beta<\alpha$ the set
${\Cal V}^{\beta}_1 = \{U^{\beta}_n: n\in S_{\beta} \mathrm{\ and\ } x_1\in U^{\beta}_n\}$ is an
$\omega$-fat family. For each $n$ choose $x_{n+1}\in J_{n+1}\cap H_{\alpha}\setminus\{x_1,\dots,x_n\}$
such that ${\Cal V}^{\beta}_{n+1}:= \{U^{\beta}_m\in{\Cal V}^{\beta}_n: x_{n+1}\in U^{\beta}_m\}$ is an
$\omega$-fat family. Finally apply Lemma \ref{ppointlemma} to choose for each $\beta<\alpha$ an
$\omega$-fat family ${\Cal V}_{\beta}\sbst {\Cal V}^{\beta}_1$ such that for each $n$
${\Cal V}^{\beta}\sbst^*{\Cal V}^{\beta}_n$, and set $X_{\alpha} = \{x_n:n\in\naturals\}$.
Observe that each ${\Cal V}^{\beta}$ is a $\g$-cover of
$\cup_{\g\leq\alpha}X_\g\cup\rationals$, and $X_{\alpha}$ is a dense subset of $\reals$.
For each $\beta<\alpha$ define $S_{\alpha}(\beta):=\{m: U^{\beta}_m\in {\Cal V}^{\beta}\}$.}
\item[Yes:]{ Then first choose for each $n$ a $U^{\alpha}_n\in{\Cal U}^{\alpha}_n$ such that
$\{U^{\alpha}_n:n\in\naturals\}$ is a $\g$-cover of
$\cup_{\g<\alpha}X_\g\cup\rationals$. For each $\beta<\alpha$ set
$S_{\beta}(\alpha) = \naturals$. Next we construct $X_{\alpha}$. For convenience, put
$H_{\alpha} = G_{\alpha}\setminus(\cup_{\g<\alpha}X_\g\cup\rationals)$.
Applying Lemma \ref{targetinglemma} choose $x_1\in J_1\cap H_{\alpha}$ such that for each $\beta<\alpha$
the set ${\Cal U}^{\beta}_1:=\{U^{\beta}_n:n\in S_{\beta} \mathrm{\ and\ }x_1\in U^{\beta}_n\}$ is
$\omega$-fat, and ${\Cal U}^{\alpha}_1 = \{U^{\alpha}_n: x_1\in U^{\alpha}_n\}$ is $\omega$-fat.
For each $n$ choose $x_{n+1}\in J_{n+1}\cap H_{\alpha}\setminus\{x_1,\,\dots,\,x_n\}$ such that for
$\beta\leq\alpha$ we have
${\Cal V}^{\beta}_{n+1} = \{U^{\beta}_m\in{\Cal V}^{\beta}_n: x_{n+1}\in U^{\beta}_m\}$ is an
$\omega$-fat family. Finally, by Lemma \ref{ppointlemma} choose for each $\beta$ an $\omega$-fat family
${\Cal V}^{\beta}$ such that for all $n$ ${\Cal V}^{\beta}\sbst^*{\Cal V}^{\beta}_n$. Observe that each
${\Cal V}^{\beta}$ is a $\g$-cover of $\cup_{\beta\le \alpha}X_{\beta}\cup\rationals$. For
$\beta\le\alpha$ define: $S_{\alpha}(\beta) = \{n:U^{\beta}_n\in{\Cal V}^{\beta}\}$.}
\end{itemize}

In either case we succeeded in extending the satisfiability of the recursive requirements before stage $\alpha$, to stage $\alpha$.
\end{proof}

\begin{cor}[CH]\label{sonebombomandufingammagamma}
There is a set of real numbers with property $\sone(\BO,\BO)$ which does not have property
$\ufin(\Gamma,\Gamma)$.
\end{cor}
\begin{proof}
We may think of having carried out the preceding construction in $\NN$; here, every set with property $\ufin(\Gamma,\Gamma)$ is bounded, and so meager. But a Lusin set is non-meager.
\end{proof}

\subsection*{Special elements of $\sone(\BG,\BG)$}

Our next task is to determine the relationship of the top row of Figure \ref{boreldiagr} to the bottom rest of
Figure \ref{basicdiagr}. For this we compare $\sone(\BG,\BG)$ with $\sone(\O,\O)$ and with $\sfin(\Omega,\Omega)$.
A set $X$ of real numbers
is said to be a  \emph{Sierpi\'nski set} if it is uncountable, and its intersection with each Lebesgue measure
zero set is countable. More generally, for an uncountable cardinal number $\kappa$ a set of real numbers is a
$\kappa$-Sierpi\'nski set if it has cardinality at least $\kappa$, but its intersection with each set of
Lebesgue measure zero is less than $\kappa$.

In Theorem 2.9 of \cite{jmss} it was shown that all Sierpi\'nski sets have the property $\ufin(\BG,\BG)$.
This also follows easily from our characterization of $\sone(\BG,\BG)$ (Theorem \ref{bgammadualtob}),
since each countable set has this property.
Indeed, our characterization and the fact that every set of real numbers of cardinality less than $\mathfrak{b}$
has property $\sone(\BG,\BG)$ gives that every $\b$-Sierpi\'nski set has property $\sone(\BG,\BG)$.
Since sets of real numbers having property $\sone(\O,\O)$ have measure zero, no $\b$-Sierpi\'nski
set has property $\sone(\O,\O)$.

\forget
\forgotten

\forget
{\Huge
We must explain why we are working in $\R$ rather than $\NN$ as we did with the Lusin sets!!!
}{\huge: Agreed, but I think at this stage we should perhaps not make it part of this paper -- I am not convinced anymore that it is that relevant to our topic. There seems to be something subtle going on, and perhaps when we understand it - and it seems important - we could write about it.}{\emph{I propose we drop this task to bring this paper to a conclusion??}}
\forgotten

Let $\P$ denote the set of irrational numbers.

\begin{lem}\label{covnconstruction}
If $\cov({\Cal N}) = \cof({\Cal N})$, and if $Y\sbst\P$
has cardinality at most $\cof({\Cal N})$, then there is a $\cov({\Cal N})$-Sierpi\'nski set $S\sbst\P$
such that $Y\sbst S + S\sbst\P$.
\end{lem}
\begin{proof}
Let $\{y_{\alpha}:\alpha<\cov({\Cal N})\}$ enumerate $Y$. Let $\{N_{\alpha}:\alpha<\cov({\Cal N})\}$
enumerate a cofinal family of meager sets, and construct $S$ recursively as follows: At stage $\alpha$ set
$$
X_{\alpha} = \bigcup_{i<\alpha}\left(
\{a_i, b_i\}\cup (\Q-a_i) \cup (\Q-b_i)\cup N_i 
\right).
$$
Note that for each $x\in  \P\setminus X_{\alpha}$ and $i<\alpha$, $x+a_i$ and $x+b_i$ are irrational.

$X_{\alpha}$ is a union of fewer than $\cov({\Cal N})$ measure zero sets.
As in Lemma \ref{covmconstruction}, we can choose $a_\alpha,b_\alpha\in \P\setminus X_{\alpha}$ such that
$a_{\alpha} + b_{\alpha} = y_{\alpha}$. (Note that $y_\alpha\in\P$.)

Finally, set $S = \{a_{\alpha}:\alpha<\cov({\Cal N})\}\cup\{b_{\alpha}:\alpha<\cov({\Cal N})\}$.
Then $S$ is a $\cov({\Cal N})$-Sierpi\'nski set and $Y\sbst S+S\sbst\P$.
\end{proof}

\begin{thm} \label{sonebgbgandsoneoo}
If $\b = \cov({\Cal N}) = \cof({\Cal N})$, then
there is a $\b$-Sierpi\'nski set of real numbers $S$ such that:
\be
\i $S$ has property $\sone(\BG,\BG)$,
\i $S$ does not have property $\sone(\O,\O)$,
\i $S\x S$ does not have property $\ufin(\Gamma, \O)$,
\i $S$ does not have property $\sfin(\Omega, \Omega)$.
\ee
\end{thm}
\begin{proof}
Note that the hypothesis $\b = \cof({\Cal N})$ implies that $\b = \d$.
Let $\Psi$ be a homeomorphism from the irrationals onto $\NN$.
Let $D\sbst{}\NN$ be a dominating family of size $\d$, and  set $Y=\Psi^{-1}[D]$.
Use Lemma \ref{covnconstruction} to construct a $\b$-Sierpi\'nski set  $S\sbst\P$
such that $Y\sbst S+S\sbst\P$.
Now, define $f:S\x S:\to{}\NN$ by $f(x,y)=\Psi(x+y)$. Then $f$ is
continuous, and $f[S\x S] = \Psi[S+S]\spst \Psi[X] = D$ is dominating.
This makes 1,2, and 3.

Now, in \cite{jmss} it is proved that $\sfin(\Omega, \Omega)$ is closed under taking
finite powers. Thus, 4 follows from 3.
\end{proof}

Thus, we have that $\sone(\BG,\BG)$ is not provably contained in $\sfin(\Omega, \Omega)$.
It follows that Figure \ref{boreldiagr} gives all the provable relations among the Borel covering classes.

In light of Theorem \ref{borelhurewicz},
the following Theorem of Rec{\l}aw \cite{irek2} implies that
none of the properties involving open classes implies any of the properties involving
Borel classes.
Rec{\l}aw's proof assumes Martin's axiom, but the partial order used is
$\s$-centered so that in fact $\p = \c$ is enough.

\begin{thm}[$\p=\c$]\label{gamma-onto-reals}
There is a set having the $\sone(\Omega,\Gamma)$ property 
which can be mapped onto $\NN$ by a Borel function.
\end{thm}

\begin{figure}
\unitlength=.30mm
\begin{picture}(300.00,500.00)(0,0)
\put(53.00,55.00){\makebox(0,0)[cc]
{{\footnotesize $\sone({\Omega},{\Gamma})$ }}}
\put(20.00,-5.00){\makebox(0,0)[cc]
{{\footnotesize $\sone(\BO,\BG)$ }}}
\put(53.00,235.00){\makebox(0,0)[cc]
{{\footnotesize $\sone({\Gamma},{\Gamma})$ }}}
\put(20.00,175.00){\makebox(0,0)[cc]
{{\footnotesize $\sone(\BG,\BG)$ }}}
\put(119.00,355.00){\makebox(0,0)[cc]
{{\footnotesize $\ufin({\Gamma},{\Gamma})$} }}
\put(50.00,-3.00){\vector(1,0){100.00}}
\put(50.00,176.0){\vector(1,0){100.00}}
\put(210.0,-3.00){\vector(1,0){100.00}}
\put(210.0,176.0){\vector(1,0){100.00}}
\put(83.00,237.0){\vector(1,0){100.00}}
\put(243.0,237.0){\vector(1,0){100.00}}
\put(149.0,355.0){\vector(1,0){100.00}}
\put(309.0,355.0){\vector(1,0){100.00}}
\put(20.00,5.00){\vector(0,1){150}}
\put(340.0,5.00){\vector(0,1){150}}
\put(373.0,65.0){\vector(0,1){150}}
\put(180.0,5.00){\vector(0,1){70}}
\put(180.0,93.0){\vector(0,1){70}}
\put(25.00,3.00){\vector(1,2){20.00}}
\put(185.0,3.00){\vector(1,2){20.00}}
\put(345.0,3.00){\vector(1,2){20.00}}
\put(25.00,183.0){\vector(1,2){20.0}}
\put(185.0,183.0){\vector(1,2){20.0}}
\put(345.0,183.0){\vector(1,2){20.0}}
\put(218.0,63.00){\vector(1,2){20.0}}
\put(218.0,243.0){\vector(1,2){20.0}}
\put(251.0,303.0){\vector(1,2){20.0}}
\put(58.00,243.0){\vector(1,2){50.0}}
\put(378.0,243.0){\vector(1,2){50.0}}
\put(53.00,60.00){\line(0,1){110.0}}
\put(53.00,180.0){\vector(0,1){40.0}}
\put(213.0,60.00){\line(0,1){110.0}}
\put(213.00,180.0){\vector(0,1){40.0}}
\put(246.0,125.00){\line(0,1){43.0}}
\put(246.0,185.00){\line(0,1){43.0}}
\put(246.0,245.00){\vector(0,1){41.0}}
\put(70.00,55.00){\line(1,0){100}}
\put(183.0,55.00){\vector(1,0){10}}
\put(230.00,55.00){\line(1,0){100}}
\put(343.0,55.00){\vector(1,0){10}}
\put(213.00,55.00){\makebox(0,0)[cc]
{{\footnotesize $\sone({\Omega},{\Omega)}$ }}}
\put(180.00,-5.00){\makebox(0,0)[cc]
{{\footnotesize $\sone(\BO,\BO)$ }}}
\put(213.00,235.00){\makebox(0,0)[cc]
{{\footnotesize $\sone({\Gamma},{\Omega})$ }}}
\put(180.00,175.00){\makebox(0,0)[cc]
{{\footnotesize $\sone(\BG,\BO)$ }}}
\put(246.00,295.00){\makebox(0,0)[cc]
{{\footnotesize $\sfin({\Gamma},{\Omega})$ }}}
\put(279.00,355.00){\makebox(0,0)[cc]
{{\footnotesize $\ufin({\Gamma},{\Omega})$ }}}
\put(246.00,115.00){\makebox(0,0)[cc]
{{\footnotesize $\sfin({\Omega},{\Omega})$ }}}
\put(180.00,85.00){\makebox(0,0)[cc]
{{\footnotesize $\sfin(\BO,\BO)$ }}}
\put(373.00,55.00){\makebox(0,0)[cc]
{{\footnotesize $\sone(\O,\O)$ }}}
\put(340.00,-5.00){\makebox(0,0)[cc]
{{\footnotesize $\sone(\B,\B)$ }}}
\put(373.00,235.00){\makebox(0,0)[cc]
{{\footnotesize $\sone({\Gamma},\O)$ }}}
\put(340.00,175.00){\makebox(0,0)[cc]
{{\footnotesize $\sone(\BG,\B)$ }}}
\put(439.00,355.00){\makebox(0,0)[cc]
{{\footnotesize $\ufin({\Gamma},\O)$ }}}
\end{picture}
\vspace{0.25in}

\caption{The Combined Diagram \label{combineddiagr2}}
\end{figure}
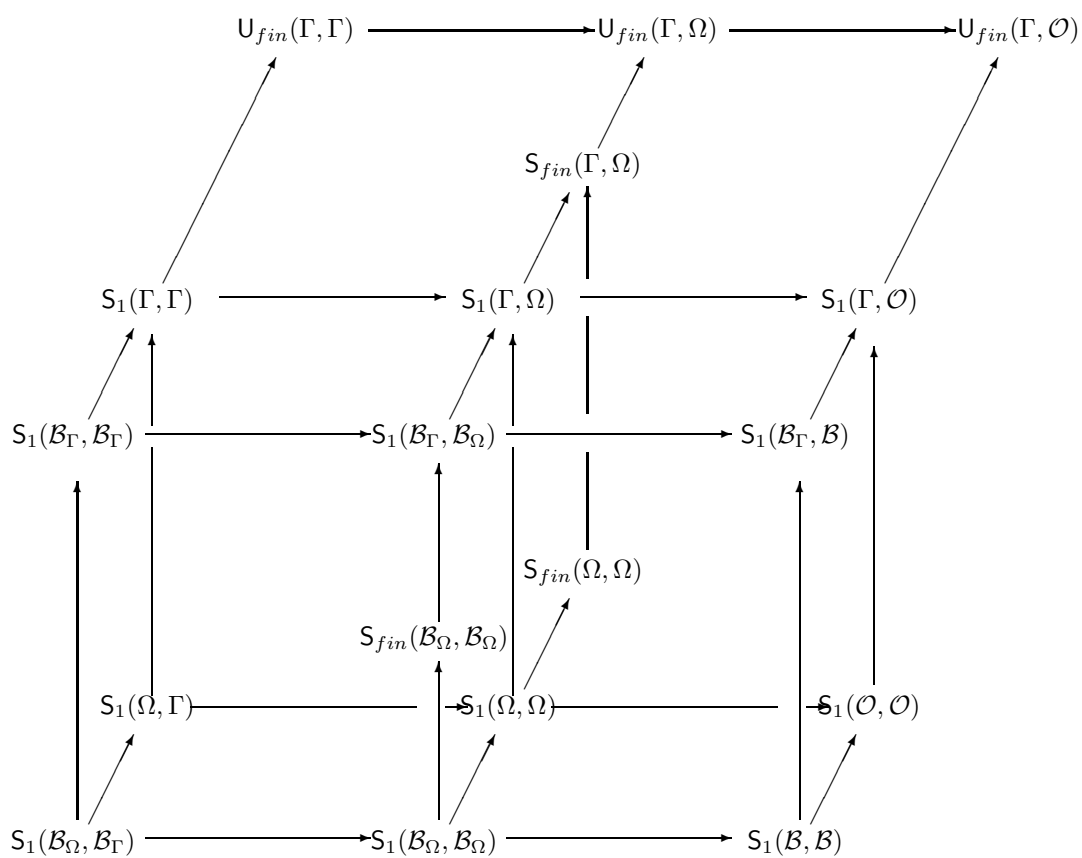

Figure \ref{combineddiagr2} summarizes the relationships among the various classes considered so far in this paper and
in \cite{jmss}, including the Borel classes. In this diagram there must also be a vector pointing from $\sfin(\BO,\BO)$ to
$\sfin(\Omega,\Omega)$; we omitted this one for ``aesthetic" reasons.

With this we have now shown that in Figure \ref{combineddiagr2}, no arrows can be added
to, or removed from, the layer of Borel classes.



At present it is not known if there always is an uncountable set of real numbers which belongs
to some class in Figure \ref{boreldiagr}.
In light of what we know about this diagram, the most modest form of this question is
\begin{prob}
Is there always an uncountable set of reals with property $\sone(\BG,\B)$?
\end{prob}
while the boldest form would be
\begin{prob}
Is there always an uncountable set of real numbers with property $\sfin(\BO,\BO)$?
\end{prob}

\subsection*{Special elements of $\sone(\BO,\BG)$}

It might be wondered whether any of our Borel notions trivializes to contain only
sets of size smaller than the critical cardinality of that notion.
With the knowledge obtained thus far, the only candidate to trivialize is
$\sone(\BO,\BG)$. A Theorem of Brendle \cite{JORG} shows that this is not the case.

\begin{thm}[CH]
There is a set of reals $X$ of size $\c(=\aleph_1)$ which has property $\sone(\BO,\BG)$.
\end{thm}


\section{Preservation of properties}

The selection properties for open covers are preserved when taking
continuous images or closed subsets \cite{jmss}.
We have the following analogue.

\begin{thm}\label{borelimages}  Let ${\Pi}$ be one of
 $\sone$, $\sfin$, or $\ufin$ and let $\UU$ and
 $\V$ range over the set $\{\B,\BO,\B_\Lambda,\BG\}$.
Assume that $X$ has property ${\Pi}(\UU,\V)$. Then:
\be
\i If $Y$ is a Borel subset of $X$, then $Y$ has property  ${\Pi}(\UU,\V)$;
\i If $f:X\to Y$ is Borel and onto, then $Y$ has property  ${\Pi}(\UU,\V)$.
\ee
\end{thm}
\begin{proof}
This proof is similar to the proof of Theorem 3.1 in \cite{jmss}.
\end{proof}

In particular, if $\UU$ and $\V$ are among $\{\O, \Omega, \Lambda, \Gamma\}$ for $X$, and $X$ has property
$\Pi(\B_\UU,\B_\V)$ for some $\Pi$, then every Borel image of $X$ has property $\Pi(\UU,\V)$. This gives rise
to the following question:
Using the above notation, assume that every Borel image of $X$ has property $\Pi(\UU,\V)$. Does $X$ necessarily
have the $\Pi(\B_\UU,\B_\V)$ property?
For the following classes, a positive answer was given:
\begin{itemize}
\i $\sone(\O,\O)$ -- Theorem \ref{borelimagec''}.
\i $\ufin(\Gamma,\Gamma)$ -- Theorem \ref{hure1}.
\i $\sone(\Gamma,\Gamma)$ -- this one follows from the preceding one, since $\sone(\Gamma,\Gamma)$
implies $\ufin(\Gamma,\Gamma)$, and $\sone(\BG,\BG)$ is equivalent to $\ufin(\BG,\BG)$
(Theorem \ref{borelbounded}).
\i $\ufin(\Gamma,\O)$ -- Theorem \ref{hure2}.
\i $\sone(\Gamma,\O)$ -- this one too follows from the preceding one, since $\sone(\Gamma,\O)$
implies $\ufin(\Gamma,\O)$, and $\sone(\BG,\B)$ is equivalent to $\ufin(\BG,\B)$ (Theorem \ref{borelhurewicz}).
\i $\sone(\Omega,\Gamma)$ -- Theorem \ref{omegagammareclaw}.
\end{itemize}
For the following classes, the problem remains open:
\begin{itemize}
\i $\sone(\Gamma,\Omega)$, $\sfin(\Gamma,\Omega)$, and $\ufin(\Gamma,\Omega)$ --
If 4 implies 3 were true in Remark \ref{openmaxfin}, we could have added these classes to the positive list.
\i $\sone(\Omega,\Omega)$.
\i $\sfin(\Omega,\Omega)$.
\end{itemize}

\subsection*{Finite powers}


$\sone(\B,\B)$ is not provably closed under taking finite powers.
\begin{thm}\label{sonebbandpowers}
If $\cov({\Cal M}) = \cof({\Cal M})$, then there exists a set of reals $X$ such that
$X$ has property $\sone(\B,\B)$, and $X\x X$ does not have the property $\ufin(\Gamma,\O)$.
\end{thm}
\begin{proof}
The $\cov(\M)$-Lusin set $L$ from Theorem \ref{sonebbnotufingamom} has the property that
$L+L$, a continuous image of $L\x L$, is dominating.
Thus, $L\x L$ does not have the property $\ufin(\Gamma,\O)$.
\end{proof}

Dually, Theorem \ref{sonebgbgandsoneoo} shows that $\sone(\BG,\BG)$
is not provably closed under taking finite powers.

\begin{prob}
Is any of the classes $\sone(\BO,\BG)$, $\sone(\BO, \BO)$, and
$\sfin(\BO, \BO)$ closed under taking finite powers?
\end{prob}

Note that a positive answer to Problem \ref{sonebombompowersareinbb}
would imply that $\sone(\BO,\BO)$ is closed under taking finite powers.
Similarly, a positive answer to Problem \ref{sfinbombompowersareinbgb}
would imply that $\sfin(\BO,\BO)$ is closed under taking finite powers.

\section{Connections with other approaches to smallness properties}\label{connections}

   Three schemas for describing smallness of sets of real numbers have been
   developed over recent years. These have their roots in classical
   literature and can be described, broadly speaking, by:
\begin{itemize}
\item{Properties of the vertical sections of a
   sufficiently describable planar set;}
\item{Properties of the image in $\NN$ under a
   sufficiently describable function;}
\item{Selection properties for sequences of sufficiently describable
   topologically significant families of subsets.}
\end{itemize}

   The vertical sections schema has been inspired by the papers
\cite{pawlikowski}, \cite{raisonnier} and \cite{irek}, and is as follows:

   Let $H$ be a subset of $\reals\times\reals$ and let $\J$ be a collection
   of subsets of $\reals$. For $x$ and $y$ real numbers, define
\[H_x = \{y\in\reals: (x,y)\in H\};
\]
\[H^y = \{x\in\reals: (x,y)\in H\}.
\]
   A Borel set $H$ is said to be a $\J$-set if for each $x$ $H_x\in \J$.

   The following three collections of subsets of the real line have been
   defined in terms of properties of vertical sections -- see
   \cite{pawlikowskireclaw}:
\begin{itemize}
\item[${\sf ADD}(\J)$:]
{The set of $X\sbst\reals$ such that for each $\J$-set $H$, $\cup_{x\in X}H_x\in \J$;}
\item[${\sf COV}(\J):$]
{The set of $X\sbst \reals$ such that for each $\J$-set $H$, $\cup_{x\in X}H_x\neq\reals$;}
\item[${\sf COF}(\J):$]
{The set of $X\sbst\reals$ such that $\{H_x:x\in X\}$ is not a cofinal subset of $\J$.}
\end{itemize}

The sets in ${\sf COV}({\Cal M})$ have also been called $R^{\Cal
M}$-sets in \cite{BarJu}; in that paper it was shown
that $X$ is an $R^{\Cal M}$-set if, and only if, every Borel image
of $X$ in $\NN$ has property $\sone(\O,\O)$. It was shown in
\cite{BarSch} that this class is also characterized by
$\sone(\B,\B)$.

The sets in ${\sf ADD}({\Cal M})$ have also been called $SR^{\Cal
M}$-sets, and it has been shown in \cite{BarJu} that
$X$ is in ${\sf ADD}({\Cal M})$  if, and only if, every Borel
image of $X$ in $\NN$ has both properties $\sone(\O,\O)$ and
$\ufin(\Gamma,\Gamma)$. Due to a result in \cite{nsw}, a set $X$
of real numbers has both properties $\sone(\O,\O)$ and
$\ufin(\Gamma,\Gamma)$ if, and only if, it has the property $(*)$
which was introduced in \cite{gerlitsnagy}. Using our results here
and results of \cite{nsw} one can show that a set of reals has
property ${\sf ADD}({\Cal M})$ if, and only if,
it is a member of $\sone(\B,\B)$ and $\sone(\BG,\BG)$.

    The ``properties of the image" schema takes inspiration from three papers \cite{hurewicz27}, \cite{irek} and
\cite{rothberger41} (Lemma 3). In each of these papers it is proven that a set of real numbers has a certain property
of interest if, and only if, each of its continuous images (in some cases into a specific range space) has another property
of interest.

The following four classes of sets were introduced in \cite{pawlikowskireclaw}:
\begin{itemize}
\item[${\sf NON}(\J)$:]{ The set of $X\sbst\reals$ such that for every Borel
   function $f$ from $\reals$ to $\reals$, $f[X]$ is a member of $\J$;}
\item[{\sf P}:]{ The set of $X\sbst\reals$ such that for no Borel
   function $f$ from $\reals$ to $\inf$, $f[X]$ is a power;}
\item[{\sf B}:]{ The set of $X\sbst\reals$ such that for every Borel
   function $f$ from $\reals$ to $\NN$, $f[X]$ is bounded
   under eventual domination;}
\item[{\sf D}:]{ The set of $X\sbst \reals$ such that for every Borel
   function $f$ from $\reals$ to $\NN$, $f[X]$ is not a
   dominating family.}
\end{itemize}

\begin{figure}
\unitlength=.95mm
\begin{picture}(140.00,60.00)(10,10)
\put(20.00,20.00){\makebox(0,0)[cc]
{${\sf ADD}({\Cal N})$ }}
\put(60.00,20.00){\makebox(0,0)[cc]
{${\sf ADD}({\Cal M})$ }}
\put(100.00,20.00){\makebox(0,0)[cc]
{${\sf COV}({\Cal M})$ }}
\put(140.00,20.00){\makebox(0,0)[cc]
{${\sf NON}({\Cal N})$ }}
\put(20.00,60.00){\makebox(0,0)[cc]
{${\sf COV}({\Cal N})$ }}
\put(60.00,60.00){\makebox(0,0)[cc]
{${\sf NON}({\Cal M})$ }}
\put(100.00,60.00){\makebox(0,0)[cc]
{${\sf COF}({\Cal M})$ }}
\put(140.00,60.00){\makebox(0,0)[cc]
{${\sf COF}({\Cal N})$ }}
\put(100.00,40.00){\makebox(0,0)[cc]
{${\sf D}$ }}
\put(60.00,40.00){\makebox(0,0)[cc]
{${\sf B}$ }}
\put(28.00,61.00){\vector(1,0){20.00}}
\put(110.00,20.00){\vector(1,0){20.00}}
\put(110.00,61.00){\vector(1,0){20.00}}
\put(70.00,61.00){\vector(1,0){20.00}}
\put(28.00,20.00){\vector(1,0){20.00}}
\put(70.00,20.00){\vector(1,0){20.00}}
\put(70.00,40.00){\vector(1,0){20.00}}

\put(20.00,25.00){\vector(0,1){29.00}}
\put(60.00,25.00){\vector(0,1){11.00}}
\put(60.00,42.00){\vector(0,1){12.00}}
\put(100.00,25.00){\vector(0,1){11.00}}
\put(100.00,42.00){\vector(0,1){12.00}}
\put(140.00,25.00){\vector(0,1){29.00}}
\end{picture}
\caption{Chichon-like Diagram \label{chichondiagr}}
\end{figure}
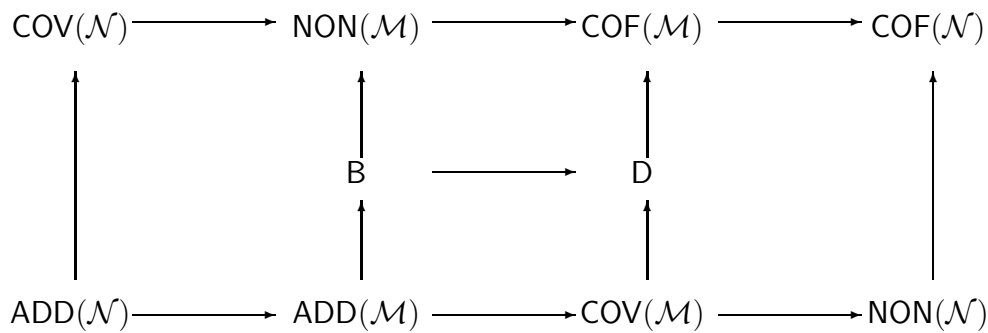

The classes of sets defined by these two schemas are related for the special case
where $\J$ is ${\Cal M}$, the collection of
meager sets of real numbers, or ${\Cal N}$, the collection of measure zero
subsets of the real line. The results from
\cite{pawlikowskireclaw} regarding the interrelationships of these classes of sets
are summarized in Figure \ref{chichondiagr}.

The relationship between Figure \ref{chichondiagr} and the well-known Chichon diagram that expresses provable
relationships among certain cardinal numbers is that a cardinal number in a particular position in Cichon's
diagram is actually the minimal cardinality for a set of real numbers not belonging to the class in the
corresponding position in Figure \ref{chichondiagr}.

Our results imply the following.
\begin{cor}
$\COF(\M)$ contains a set of reals whose size is $\cov(\M)$.
\end{cor}
\begin{proof}
If $\cov(\M)<\cof(\M)(\ =\unif(\COF(\M))\ )$, then any set of size $\cov(\M)$ will do.
Otherwise by Theorem \ref{sonebbnotufingamom} there exists a $\cov(\M)$-Lusin set in $\sone(\B,\B)$,
which is in $\COV(\M)$.
\end{proof}

In \cite{hurewicz27} Hurewicz characterized the covering properties $\ufin(\Gamma,\Gamma)$ and $\sfin(\O,\O)$ in terms of properties of the continuous images in $\NN$. In particular, Hurewicz showed that $X$ has property $\ufin(\Gamma,\Gamma)$ if, and only if, each continuous image of $X$ in $\NN$ is bounded. He also showed that $X$ has property $\sfin(\O,\O)$ if, and only if, each continuous image of $X$ into $\NN$ is not a dominating family.
The sets in {\sf B} have also been called \emph{A-sets} in
\cite{BarSch}; where they show that that ${\sf B} =
\ufin(\BG,\BG)$, and ${\sf D} = \sfin(\B,\B)$. By our results here
we know ${\sf B} = \sone(\BG,\BG)$, and ${\sf D} = \sone(\BG,\B)$.

\end{document}